\documentclass[nohyperref]{article}
\usepackage{microtype}
\usepackage{graphicx}
\usepackage{booktabs} 

\usepackage{hyperref}

\usepackage[accepted]{myarxiv}

\usepackage{amsmath}
\usepackage{amssymb}
\usepackage{mathtools}
\usepackage{amsthm}
\usepackage{pifont}

\usepackage[utf8]{inputenc} 
\usepackage[T1]{fontenc}    
\usepackage{hyperref}       
\usepackage{url}            
\usepackage{booktabs}       
\usepackage{amsfonts}       
\usepackage{nicefrac}       
\usepackage{microtype}      
\usepackage{xcolor}         

\usepackage{multicol}
\usepackage{microtype}
\usepackage{graphicx}
\usepackage{amsmath}
\usepackage{amssymb}
\usepackage{mathtools}
\usepackage{amsthm}
\usepackage{subfigure}

\usepackage[textsize=tiny]{todonotes}
\usepackage{epsfig}
\usepackage{cancel,soul,ulem}
\usepackage{algorithm}
\usepackage{algorithmic}

\theoremstyle{plain}
\newtheorem{theorem}{Theorem}[section]
\newtheorem{proposition}[theorem]{Proposition}
\newtheorem{lemma}[theorem]{Lemma}

\theoremstyle{definition}
\newtheorem{definition}[theorem]{Definition}
\newtheorem{assumption}[theorem]{Assumption}
\theoremstyle{remark}

\newcommand{\la}{\langle}
\newcommand{\ra}{\rangle}
\newcommand{\noi}{\noindent}
\newcommand{\nn}{\nonumber}
\newcommand{\beq}{\begin{eqnarray}}
\newcommand{\eeq}{\end{eqnarray}}

\def\etas{\boldsymbol{\eta}}

\def\x{\mathbf{x}}
\def\c{\mathbf{c}}
\def\y{\mathbf{y}}
\def\z{\mathbf{z}}
\def\g{\mathbf{g}}
\def\d{\mathbf{d}}

\def\w{\mathbf{w}}
\def\v{\mathbf{v}}

\def\G{\mathbf{G}}
\def\A{\mathbf{A}}
\def\I{\mathbf{I}}

\def\QQ{\mathbf{Q}}

\def\E{\mathbb{E}}

\def\R{\mathcal{R}}
\def\Q{\mathcal{Q}}
\def\P{\mathcal{P}}
\def\K{\mathcal{K}}

\def\J{\mathcal{J}}
\def\U{\mathcal{U}}

\def\ei{e_i}
\def\eit{e_{i^t}}

\def\red{\textcolor{red}}
\def\conenn{\textcolor[rgb]{0.9961,0,0}}

\newcommand{\cone}[1]{\textbf{\conenn{#1}}}

\newcommand{\bfit}[1]{\textit{\textbf{#1}}}

\usepackage[capitalize,noabbrev]{cleveref}

\icmltitlerunning{Coordinate Descent Methods for Fractional Minimization }

\begin{document}

\twocolumn[

\icmltitle{Coordinate Descent Methods for Fractional Minimization}



\icmlsetsymbol{equal}{*}

\begin{icmlauthorlist}
\icmlauthor{Ganzhao Yuan}{equal,yyy}
\end{icmlauthorlist}

\icmlaffiliation{yyy}{Peng Cheng Laboratory, China}

\icmlcorrespondingauthor{Ganzhao Yuan}{yuangzh@pcl.ac.cn}

\icmlkeywords{Machine Learning, ICML}

\vskip 0.3in
]



\printAffiliationsAndNotice{ } 

\begin{abstract}
We consider a class of structured fractional minimization problems, in which the numerator part of the objective is the sum of a differentiable convex function and a convex non-smooth function, while the denominator part is a convex or concave function. This problem is difficult to solve since it is non-convex. By exploiting the structure of the problem, we propose two Coordinate Descent (CD) methods for solving this problem. The proposed methods iteratively solve a one-dimensional subproblem \textit{globally}, and they are guaranteed to converge to coordinate-wise stationary points. In the case of a convex denominator, under a weak \textit{locally bounded non-convexity condition}, we prove that the optimality of coordinate-wise stationary point is stronger than that of the standard critical point and directional point. Under additional suitable conditions, CD methods converge Q-linearly to coordinate-wise stationary points. In the case of a concave denominator, we show that any critical point is a global minimum, and CD methods converge to the global minimum with a sublinear convergence rate. We demonstrate the applicability of the proposed methods to some machine learning and signal processing models. Our experiments on real-world data have shown that our method significantly and consistently outperforms existing methods in terms of accuracy.
\end{abstract}

\section{Introduction}

Fractional optimization, referring to the problem of minimizing or maximizing an objective involving one or more ratios of functions, has been extensively studied for decades. Fractional optimization problem is widely used in machine learning, signal processing, economics, wireless communication and many other fields. Three classes of factional optimization problems for minimizing the ratio of two functions are extensively investigated in the literature. They are named according to the functions in the numerator and denominator: \textit{\textbf{(i)}} linear fractional problems if both functions are linear; \textit{\textbf{(ii)}} convex-convex fractional problems if both functions are convex; \textit{\textbf{(iii)}} convex-concave fractional problems if the numerator is convex and the denominator is concave. We refer the readers to \cite{stancu2012fractional,schaible1995fractional} for an overview.

This paper mainly focuses on the following convex-convex or convex-concave Fractional Minimization Problem (FMP) (`$\triangleq$' means define):
\beq \label{eq:main}
 \bar{\x} \in \arg \min_{\x\in\mathbb{R}^n}~ F(\x)\triangleq \frac{f(\x) + h(\x)}{g(\x)}  .
\eeq

\noi We impose the following assumptions on Problem (\ref{eq:main}) throughout this paper. \textit{\textbf{(A-i)}} $F(\x)$ only takes finite values, and it always holds that: $f(\x) + h(\x)\geq0$ and $g(\x)>0$ for all $\x$. \textit{\textbf{(A-ii)}} $f(\cdot)$ is convex and differentiable, and its gradient is coordinate-wise Lipschitz continuous with constant $\c_i\geq0$ that \cite{nesterov2012efficiency,nesterov2003introductory}:
\beq \label{eq:f:Lipschitz}
\textstyle f(\x + \eta \ei) \leq \Q_i(\x,\eta) \triangleq f(\x) +  \nabla_i f(\x) \eta + \frac{\c_i}{2}\eta^2
\eeq
\noi $\forall \x,\eta, i=1,...,n$. Here $\c \in \mathbb{R}^n$, and $\ei \in \mathbb{R}^n$ is an indicator vector with one on the $i$-th entry and zero everywhere else. \textit{\textbf{(A-iii)}} $h(\cdot)$ is convex and coordinate-wise separable with $h(\x) = \sum_{i=1}^n h_i(\x_i)$. Typical examples of $h(\x)$ are the $\ell_1$ norm function and the bound constrained function. \textit{\textbf{(A-iv)}} The denominator $g(\cdot)$ is either \textit{\textbf{(i)}} a convex but not necessarily differentiable function or \textit{\textbf{(ii)}} a concave and differentiable function. Furthermore, $g(\cdot)$ has some special structure such that one of the following one-dimensional subproblems can be solved exactly and efficiently: 
\begin{align}
&\textstyle \min_{\eta} \frac{ \varrho + \rho  \eta + \frac{\c_i}{2}\eta^2 + h_i(\x+\eta \ei )}{g(\x+\eta \ei )} \nn\\
&\textstyle \min_{\eta} \varrho +  \rho \eta + \frac{\c_i}{2}\eta^2 + h_i(\x+\eta \ei ) - \gamma g(\x+\eta \ei ) \nn
\end{align}
 \noi for any $\varrho\in\mathbb{R}$, $\rho\in\mathbb{R}$, and $\gamma\in\mathbb{R}$. Problem (\ref{eq:main}) captures a variety of applications of interest, e.g., the sparse recovery problem \cite{LIHarmonic,LISIOPT2022}, the independent component analysis \cite{hyvarinen1997fast}, the $\ell_p$ norm eigenvalue problem, the regularized total least squares problem \cite{beck2006finding,amaral2005connections}, and the transmit beamforming \cite{sidiropoulos2006transmit}.

Coordinate Descent (CD) is an iterative algorithm that successively performs minimization along coordinate directions. Due to its simplicity and efficiency, it has been used for many years on the structured high dimensional machine learning and data mining applications including support vector machines \cite{hsieh2008dual}, non-negative matrix factorization \cite{hsieh2011fast}, LASSO \cite{tseng2009coordinate}. Its iteration complexity for convex problems has been well-studied \cite{nesterov2012efficiency,lu2015complexity}. Recently, its popularity continues to grow due to its strong optimality guarantees and superior empirical performance when it is applied to solve non-convex problems, including compressed sensing \cite{BeckE13,yuan2020block}, eigenvalue complementarity problem \cite{patrascu2015efficient}, DC minimization problem \cite{yuan2021coordinate}, k-means clustering \cite{nie2021coordinate}, sparse phase retrieval \cite{shechtman2014gespar}. To the best of our knowledge, this is the first time to apply CD methods for solving FMPs and study their theoretical properties and empirical behaviors.

\textbf{Contributions.} The contributions of this paper are as follows: \textit{\textbf{(i)}} We propose two CD methods for solving FMPs. The proposed methods iteratively solve a one-dimensional subproblem \textit{globally} until convergence. See Section \ref{sect:alg}. \textit{\textbf{(ii)}} For convex-convex FMPs, we prove that under suitable conditions the proposed CD methods find stronger coordinate-wise stationary points than existing methods, and they converge linearly. For convex-concave FMPs, we prove that CD methods converge to the global optimal solutions with a sublinear convergence rate. See Section \ref{sect:theory}. \textit{\textbf{(iii)}} We demonstrate the applicability of the CD methods to the applications of sparse recovery and $\ell_p$ norm eigenvalue problem. We show that the exact minimizer of each coordinate can be obtained by using an elaborate breakpoint searching procedure. Our experiments on real-world data have shown that our methods significantly and consistently outperforms existing approaches in terms of accuracy. See Section \ref{sect:exp}.


\textbf{Notations.} We use boldface lowercase letters and boldface uppercase letters to denote vectors and matrices, respectively. The Euclidean inner product between $\x$ and $\y$ is denoted by $\la \x,\y \ra$ or $\x^T\y$. We define $\|\x\| = \|\x\|_2 = \sqrt{\la \x,\x\ra}$. $\x_i$ is the $i$-th element of the vector $\x$. We define $\|\d\|_{\c}^2 \triangleq\sum_{i} \c_i\d_i^2$. $\I$ is the identity matrix of suitable size. $\text{dist}(\Omega,\Omega') \triangleq \inf_{\v \in \Omega,\v' \in \Omega'}\| \v - \v'\|$ denotes the distance between two sets.


\section{Applications}
A wide range of machine learning and signal processing models can be formulated as Problem (\ref{eq:main}). We briefly review two instances as follows.

\noi$\bullet$ \textbf{Application I: Sparse Recovery} \cite{LIHarmonic,LISIOPT2022,GotohTT18,bi2014exact}. It is a signal processing technique, which can effectively acquire and reconstruct the signal by finding the solution of the underdetermined linear system. Given a design matrix $\G \in \mathbb{R}^{m\times n}$ and an observation vector $\y\in \mathbb{R}^{m}$, sparse recovery can be formulated as the following FMP \cite{LIHarmonic}:
\beq \label{eq:app:sparse}
\textstyle \min_{\x}~\frac{\frac{1}{2}\|\G\x - \y\|_2^2 + \gamma \|\x\|_1}{ \gamma \sum_{j=1}^k |\x_{[j]}|},s.t.\|\x\|_{\infty}\leq \vartheta,
\eeq
\noi where $\x_{[i]}$ is the $i$-th largest component of $\x$ in magnitude, and $\gamma>0,\vartheta>0$ are given parameters.




\noi $\bullet$ \textbf{Application II: $\ell_p$ Norm Eigenvalue Problem}. Given arbitrary data matrices $\G \in \mathbb{R}^{m\times n}$ and $\mathbf{Q} \in \mathbb{R}^{n\times n}$ with $\mathbf{Q} \succ \mathbf{0}$, it aims at solving the following problem:
\beq \label{eq:pca0}
\bar{\v} = \arg \max_{\v} \|\G\v\|_p,~s.t.~\v^T\mathbf{Q}\v = 1
\eeq
\noi with $p\geq 1$. When $p=4$ and $\mathbf{Q}=\I$, Problem (\ref{eq:pca0}) reduces to the Independent Component Analysis (ICA) \cite{hyvarinen2000independent,zhai2020complete}; when $p=1$ and $\mathbf{Q} = \I$, Problem (\ref{eq:pca0}) is the $\ell_1$ PCA problem \cite{kim2019simple}. We have the following equivalent unconstrained FMPs:
\beq
\textstyle \bar{\x} = \arg \min_{\x}~\frac{\x^T\mathbf{Q}\x+\gamma_1}{\|\G\x\|_p+\gamma_2}, \label{eq:app:eig:1} \\
\textstyle \text{or}~\bar{\x} = \arg \min_{\x}~\frac{\x^T\mathbf{Q}\x +\gamma_3}{\|\G\x\|^2_p + \gamma_4}.\label{eq:app:eig:2}
\eeq
\noi Here, $\gamma_1$, $\gamma_2$, $\gamma_3$, $\gamma_4$ can be any nonnegative constant. The optimal solution to Problem (\ref{eq:pca0}) can be computed as $\bar{\v}=\pm \bar{\x} \cdot (\bar{\x}^T\mathbf{Q}\bar{\x})^{-\frac{1}{2}}$. Refer to Section \ref{app:sect:disc:ref} in the \bfit{Appendix} for detailed discussions.

\section{Related Work} \label{sect:related}

We present some related fractional optimization / minimization algorithms.

\noi \bfit{(i)} Dinkelbach's Parametric Algorithm (\bfit{DPA}) \cite{dinkelbach1967nonlinear} is one of the classical approaches for fractional optimization, which deals with Problem (\ref{eq:main}) by solving its associated parametric problem. By this approach, Problem (\ref{eq:main}) has an optimal solution $\x\in\mathbb{R}^n$ if and only if $\x$ is an optimal solution to the following problem: $\min_{\x}f(\x) + h(\x) - \bar{\lambda} g(\x)$, where $\bar{\lambda} = \frac{f(\bar{\x})+h(\bar{\x})}{g(\bar{\x})}$. However, the optimal objective value $\bar{\lambda}$ is unknown in general. Iterative procedures are considered to remedy this issue. \bfit{DPA} generates a sequence $\{\x^t\}$ as: $\x^{t+1} = \arg\min_{\x}~f(\x)+ h(\x)-\lambda^t g(\x)$, where $\lambda^t$ is renewed via $\lambda^t = \frac{f(\x^t)+h(\x^t)}{g(\x^t)}$. Note that the computational cost of solving the subproblem could be expensive since it is non-convex in general.

\bfit{(ii)} Proximal Gradient Algorithm (\bfit{PGA}) \cite{Radu2017} has been proposed for a similar class of fractional optimization problems where the denominator $g(\x^t)$ is differentiable, and can be suitably applied to Problem (\ref{eq:main}). The resulting algorithm generates a sequence $\{\x^t\}$ as: $\x^{t+1} = \arg \min_{\x} f(\x) + h(\x)  - \lambda^t \la \nabla g(\x^t),  \x - \x^{t}  \ra + \frac{1}{2 \eta^t} \|\x - \x^t\|_2^2$, where $\eta^t>0$ and $\lambda^t=F(\x^t)$.

\bfit{(iii)} Proximal Gradient-Subgradient Algorithm (\bfit{PGSA}) \cite{LISIOPT2022,LIHarmonic} assumes that $\nabla f(\cdot)$ is Lipschitz continuous with constant $L$ that: $\forall \x,\y,~f(\x) \leq \U(\x;\y) \triangleq f(\y) +  \la \nabla  f(\y) ,\x-\y\ra + \frac{L}{2}\|\y-\x\|_2^2$, and generates the new iterate using: $\x^{t+1} = \arg \min_{\x}~h(\x) + \U(\x;\x^t) - \lambda^t \la \g^t,\x-\x^t  \ra$, with $\g^t \in  \partial g (\x^t),~\lambda^t=F(\x^t)$.

\bfit{(iv)} Quadratic Transform Parametric Algorithm (\bfit{QTPA}) \cite{ShenY18,ShenY18a} introduces an additional variable $\beta\in \mathbb{R}$ and converts Problem (\ref{eq:main}) into the following equivalent variational reformulation: $\min_{\x}\frac{-g(\x)}{f(\x)+h(\x)}\Leftrightarrow\min_{\x,\beta}\beta^2(f(\x) + h(\x))-2\beta \sqrt{g(\x)}$. Similar to \bfit{DPA}, it generates a sequence $\{\x^t\}$ as: $\x^{t+1} = \arg \min_{\x}(\beta^t)^2 (f(\x) + h(\x)) - 2\beta^t \sqrt{g(\x)}$, where $\beta^t$ is renewed via $\beta^t={\sqrt{g(\x^t)}}/{(f(\x^t) + h(\x^t))}$. Note that this method is originally designed for solving multiple-ratio FMPs.

\bfit{(v)} Charnes-Cooper Transform Algorithm (\bfit{CCTA}) converts the original linear-fractional programming problem to a standard linear programming problem \cite{charnes1962programming}. Using the transformation $\y = \frac{\x}{g(\x)}$, $t=\frac{1}{g(\x)}$, one can convert Problem (\ref{eq:main}) into: $\min_{t,\y} t f(\y/t) + t h(\y/t),~s.t.~tg(\y/t) = 1$.

\bfit{(vi)} Other fractional optimization algorithms. A number of other fractional optimization algorithms have been studied in the literature. PGSA with line search is developed for possible acceleration for Problem (\ref{eq:main}) \cite{LIHarmonic}; An extrapolated proximal subgradient algorithm was proposed for solving a similar class of fractional optimization problems \cite{boct2021extrapolated}.


It is shown that any accumulation points of the sequence generated by all the algorithms above are critical points of Problem (\ref{eq:main}). 


\section{Proposed Coordinate Descent Methods} \label{sect:alg}

This section presents four variants of Coordinate Descent (CD) methods for solving the Fractional Minimization Problem (FMP) in Problem (\ref{algo:main}). 

\noi $\blacktriangleright$ \bfit{Raw Coordinate Descent}. In the $t$-th iteration of CD method, we minimize $F(\cdot)$ with respect to the $i^t$ coordinate while keeping the remaining coordinates $\{\x^t_j\}_{j\neq i^t}$ fixed. In other words, we solve the following one-dimensional subproblem: $\bar{\eta}^t \in \arg \min_{\eta \in \mathbb{R}}~\frac{f(\x^t + \eta \eit) + h(\x^t + \eta \eit) }{g(\x^t + \eta \eit)}$ and then update the solution via $\x^{t+1}  = \x^t  + \bar{\eta}^t \cdot \eit$. However, when $f(\cdot)$ and $g(\cdot)$ are complicated, this one-dimensional problem could be still difficult to solve. Some majorization or approximation techniques are needed to remedy this issue.

\noi $\blacktriangleright$ \bfit{Parametric Subgradient Coordinate Descent} is based on the parametric problem of Problem (\ref{eq:main}): $\min_{\x}f(\x)+ h(\x)-\lambda^t g(\x)$ with $\lambda^t=F(\x^t)$ is the current estimate of the objective value. When $g(\cdot)$ is convex, we have:
\beq \label{eq:lose:bound}
- g(\x+\eta\ei) \leq \mathcal{G}_{i^t}(\x^t,\eta) \triangleq -g(\x^t) - \la \partial g(\x^t),\eta  \ei \ra.
\eeq
\noi If we replace $f(\x^t + \eta \eit)$ and $-g(\x^t + \eta \eit)$ with their majorization functions $\Q_{i^t}(\x^t,\eta)$ and $\mathcal{G}_{i^t}(\x^t,\eta)$ while keep the term $h(\cdot)$ unchanged, we have:
\begin{align}
\textstyle \bar{\eta}^t &\in \textstyle \arg \min_{\eta} \Q_i(\x^t,\eta) + h(\x^t+\eta\ei) - \lambda^t \mathcal{G}_{i^t}(\x^t,\eta)\nn\\
\textstyle \x^{t+1}  &=\textstyle \x^t  + \bar{\eta}^t \cdot \eit.\nn
\end{align}
\noi However, the upper bound which only uses the subgradient of the nonconvex function $(-g(\x^t))$ in (\ref{eq:lose:bound}) could be loose, and it results in weak optimality of critical points for convex-convex FMPs \cite{LISIOPT2022,LIHarmonic}.

\noi $\blacktriangleright$  \bfit{Fractional Coordinate Descent} is rooted in the original fractional minimization function. It replaces $f(\x^t + \eta \eit)$ with its majorization (upper-bound) $\Q_{i^t}(\x^t,\eta)$ with $\Q_{i}(\x,\eta)\triangleq f(\x) +  \nabla_i f(\x) \eta + \tfrac{\c_i}{2}\eta^2$ while keeps the remaining two terms $h(\cdot)$ and $g(\cdot)$ unchanged, leading to the following iterative procedure:
\begin{align}
\textstyle \bar{\eta}^t &\in \textstyle  \arg \min_{\eta} \frac{\Q_i(\x^t,\eta) + h(\x^t+\eta\ei) }{g(\x^t + \eta \eit)}\nn\\
\textstyle \x^{t+1}  &= \textstyle \x^t  + \bar{\eta}^t \cdot \eit.\nn
\end{align}
\noi $\blacktriangleright$ \bfit{Parametric Coordinate Descent} is built upon the associated parametric problem of Problem (\ref{eq:main}). It replaces $f(\x^t + \eta \eit)$ with its majorization function $\Q_{i^t}(\x^t,\eta)$ while keeps the term $h(\cdot)$ and $g(\cdot)$ unchanged, resulting in the following updating scheme:
\begin{align}
\textstyle \bar{\eta}^t &\in  \textstyle \arg \min_{\eta} \Q_i(\x^t,\eta) + h(\x^t+\eta\ei) - \lambda^t g(\x^t + \eta \eit)\nn\\
\textstyle \x^{t+1}  &= \textstyle \x^t  + \bar{\eta}^t \cdot \eit.\nn
\end{align}
\noi $\blacktriangleright$\textbf{Choosing the Coordinate to Update}. There are mainly several strategies to decide which coordinate to update in the literature \cite{tseng2009coordinate}. \textit{\textbf{(i)}} Cyclic order rule runs all coordinates in cyclic order $1 \rightarrow 2 \rightarrow ... \rightarrow n \rightarrow 1$. \textit{\textbf{(ii)}} Random sampling rule randomly selects one coordinate to update. \textit{\textbf{(iii)}} Greedy rule picks coordinate $i^t$ such that $i^t = \arg \max_{j}|\bar{\d}^t_j|$ where $\bar{\d}^t = \arg \min_{\d}   \la \nabla f(\x^t) - F(\x^t) \partial g(\x^t), \d\ra + \frac{L}{2}\|\d\|_2^2 + h(\x^t + \d)$. Note that it has an equivalent form to the update rule of PGSA (see Section \ref{sect:related}) and $\bar{\d}^t=\mathbf{0}$ implies that $\x^t$ is a stationary point.

\begin{algorithm}[!t]
\caption{ {\bf Coordinate Descent Methods for Fractional Minimization.} }
\begin{algorithmic}
  \STATE Input: an initial feasible solution $\x^0$, $\theta>0$. Set $t=0$.
\FOR{$t = 0, 1,2,3 ... T$}
  \STATE \textit{\textbf{(S1)}} Use some strategy to find a coordinate $i^t \in \{1,...,n\}$ for the $t$-th iteration.

  \STATE \textit{\textbf{(S2)}} Define  \beq \label{eq:Ji}
\textstyle \J_i(\x,\eta)  \triangleq f(\x) +  \nabla_i f(\x) \eta + \frac{\c_i+\theta}{2}\eta^2 + h(\x+\eta\ei).\nn
  \eeq
  \noi Solve one of the following subproblems \emph{globally}:

  $\bullet$ Option \textit{\textbf{(I)}}: \bfit{FCD}.
\begin{eqnarray} \label{eq:subprob:nonconvex}
  \textstyle  \bar{\eta}^t   \in \P_{i^t}(\x^t),~ \P_{i}(\x) \triangleq \arg \min_{\eta} \frac{ \J_{i}(\x,\eta) }{g(\x + \eta \ei)}.
\end{eqnarray}

  $\bullet$ Option \textit{\textbf{(II)}}: \bfit{PCD}.
\beq \label{eq:subprob:nonconvex2}
&\textstyle\bar{\eta}^t  \in \P_{i^t}(\x^t)\\
&\textstyle\P_{i}(\x) \triangleq \arg \min_{\eta} \J_{i}(\x,\eta)  - F(\x) g(\x + \eta \ei).\nn
\eeq

\STATE \textit{\textbf{(S3)}} $\x^{t+1}  = \x^t  + \bar{\eta}^t \cdot \eit~~~ (\Leftrightarrow \x_{i^t}^{t+1}  = \x_{i^t}^t  + \bar{\eta}^t )$
   \ENDFOR
\end{algorithmic}
\label{algo:main}
\end{algorithm}

Due to the limitations of \textit{Raw Coordinate Descent} and \textit{Parametric Subgradient Coordinate Descent}, we only focus on \textit{Fractional Coordinate Descent} \bfit{(FCD)} and \textit{Parametric Coordinate Descent} \bfit{(PCD)} in the sequel. We formally present \bfit{FCD} and \bfit{PCD} in Algorithm \ref{algo:main}.

\noi \textbf{Remarks}. \bfit{(i)} Note that we increase $\frac{\c_i}{2}\eta^2$ to $\frac{\c_i+\theta}{2}\eta^2$ for the term $\J_i(\x,\eta)$. It can be viewed as appending a new proximal term $\frac{\theta}{2}\eta^2  = \tfrac{\theta}{2} \| (\x^t+\eta \eit) - \x^t \|_2^2$ to the numerator. As we will see later, the introduction of the proximal term $\tfrac{\theta}{2} \| (\x^t+\eta \eit) - \x^t \|_2^2$ is critically important for our theoretical analysis. \bfit{(ii)} Setting the derivative of the objective function with respect to $\eta$ to zero, we obtain the following \textit{necessary but not sufficient} optimality conditions for (\ref{eq:subprob:nonconvex}) and (\ref{eq:subprob:nonconvex2}), respectively:
\beq
\textstyle \alpha^t \triangleq \frac{\J_{i^t}(\x^t,\bar{\eta}^t)}{g(\x^{t+1})},~~\alpha^t  \partial_{i^t} g(\x^{t+1})  \textstyle\in   \partial \J_{i^t}(\x^t,\bar{\eta}^t), \label{eq:opt:111}\\
 \textstyle F(\x^t)\cdot  \partial_{i^t}  g(\x^{t+1}) \textstyle\in   \partial \J_{i^t}(\x^t,\bar{\eta}^t). \label{eq:opt:222}
\eeq
\noi \bfit{(iii)} Both \bfit{FCD} and \bfit{PCD} are applicable to solve both convex-convex FMPs and convex-concave FMPs. \bfit{(iv)} For convex-convex FMPs, the subproblems in (\ref{eq:subprob:nonconvex}) and (\ref{eq:subprob:nonconvex2}) are generally non-convex. However, using an elaborate breakpoint searching procedure, its exact minimizer can be obtained. This is the key insight into our CD methods. Existing methods mainly consider \textit{multiple-stage convex approximation} to handle the convex denominator term, only resulting in weak optimality of critical points \cite{LISIOPT2022,LIHarmonic,dinkelbach1967nonlinear}. Our methods directly optimize over the denominator term and globally solve a non-convex one-dimensional subproblem. Such a \textit{sequential nonconvex approximation} strategy leads to stronger optimality conditions. \bfit{(v)} In many situations, the exact minimizer of \bfit{PCD} is easier to obtained than that of \bfit{FCD} since the latter involves an objective function which is of fractional structure.



\section{Theoretical Analysis} \label{sect:theory}

We now provide some theoretical analysis of Algorithm \ref{algo:main}. We treat convex-convex FMPs and convex-concave FMPs separately. Due to space limit, all proofs are placed into the \bfit{Appendix}.

\subsection{Technical Preliminaries} \label{sebsect:technical}

We need some tools in non-smooth analysis including Fr\'{e}chet subdifferential, limiting (Fr\'{e}chet) subdifferential, and directional derivative \cite{Mordukhovich2006,Rockafellar2009,bertsekas2015convex}. For any extended real-valued (not necessarily convex) function $F: \mathbb{R}^n \rightarrow (-\infty,+\infty]$, its domain is defined by $\text{dom}(F)\triangleq \{\x\in\mathbb{R}^n: |F(\x)|<+\infty\}$. The Fr\'{e}chet subdifferential of $F$ at $\x\in\text{dom}(F)$, denoted as $\hat{\partial}F(\x)$, is defined as $\hat{\partial}{F}(\x) \triangleq \{\v\in\mathbb{R}^n: \lim \inf_{\z \rightarrow \x} \frac{ {F}(\z) - {F}(\x) - \la \v,\z-\x \ra  }{\|\z-\x\|}\geq 0 \}$. The limiting subdifferential of ${F}(\x)$ at $\x\in\text{dom}({F})$ is defined as: $\partial{F}(\x)\triangleq \{\v\in \mathbb{R}^n: \exists \x^k \rightarrow \x, {F}(\x^k)  \rightarrow {F}(\x), \v^k \in\hat{\partial}{F}(\x^k) \rightarrow \v,\forall k\}$. Note that $\hat{\partial}{F}(\x) \subseteq \partial{F}(\x)$. If $F(\cdot)$ is differentiable at $\x$, then $\hat{\partial}{F}(\x) = \partial{F}(\x) = \{\nabla F(\x)\}$ with $\nabla F(\x)$ being the gradient of $F(\cdot)$ at $\x$. When $F(\cdot)$ is convex, $\hat{\partial}{F}(\x)$ and $\partial{F}(\x)$ reduce to the classical subdifferential for convex functions, i.e., $\hat{\partial}{F}(\x) = \partial{F}(\x) = \{\v\in\mathbf{R}^n: F(\z)-F(\x)-\la\v,\z-\x \ra\geq 0,\forall \z\in\mathbb{R}^n\}$. Since $\partial h(\x)$ is coordinate-wise separable, we use $(\partial h(\x))_i$ to denote the subgradient of $h(\x)$ at $\x$ for the $i$-th component. The directional derivative of $F(\cdot)$ at $\x$ in the direction $\v$ is defined (if it exists) by $F'(\x;\v) \triangleq \lim_{t\rightarrow 0^{+}} \frac{1}{t} (F(\x + t \v) - F({\x}))$.


We present two kinds of stationary solutions for the non-convex non-differentiable FMP in (\ref{eq:main}).

\begin{definition}
(Critical Point, or \textit{C}-Point for short) A solution $\check{\x}$ is called a \textit{C}-point if \cite{LISIOPT2022}: $0 \in \partial \nabla f(\check{\x}) + \partial h(\check{\x}) - F(\check{\x}) \cdot \partial g(\check{\x})$.


\end{definition}



\begin{definition} \label{def:D:point}
(Directional Point, or \textit{D}-Point for short) A solution $\grave{\x}$ is called a \textit{D}-point if \cite{PangRA17}: $F'(\grave{\x};\y - \grave{\x} ) \geq 0,~\forall \y\in \text{dom}(F)$.

\end{definition}

\noi \textbf{Remarks}. \textit{\textbf{(i)}} The definition of \textit{C}-Point differs from the standard one $0\in \hat{\partial}F(\check{\x})$, and it holds that $\hat{\partial}F(\check{\x})= \frac{ \hat{\partial}( g(\check{\x}) (f+h) - (f(\check{\x})+h(\check{\x}))g)(\check{\x}) }{(g(\check{\x}))^2}\subseteq \partial F(\check{\x})$ \cite{LISIOPT2022}. \textit{\textbf{(ii)}} When $F(\cdot)$ is differentiable, the optimality of \textit{C}-point is equivalent to that of \textit{D}-point. \textit{\textbf{(iii)}} The expression $0 \in \partial F(\check{\x})$ is equivalent to $[\nabla f(\check{\x}) + \partial h(\check{\x})] \cap [F(\check{\x})\partial g(\check{\x})] \neq \emptyset $. \textit{\textbf{(iii)}} The function $g(\cdot)$ need not be differentiable since the sub-differential is always non-empty on convex functions. \textit{\textbf{(iv)}} All existing methods including \bfit{DPA}, \bfit{PGA}, \bfit{PGSA}, and \bfit{QTPA} as mentioned in Section \ref{sect:related} are only guaranteed to find a \textit{C}-point of Problem (\ref{eq:main}).




We make the following assumption which will be used in our theoretical analysis.

\begin{assumption}\label{ass:3}
(\bfit{Boundedness of the Denominator}) There exists a constant $\bar{g}>0$ such that $\forall \x \in \{\z~|~F(\z) \leq F(\x^0)\},~g(\x)\leq \bar{g}$.
\end{assumption}

\noi \textbf{Remarks}. As multiplying the numerator and the denominator of $F(\x)$ simultaneously by a positive constant does not change the value of $F(\x)$, Assumption \ref{ass:3} is reasonable.



We develop the following useful lemmas for both convex-convex FMPs and convex-concave FMPs.

\begin{lemma} \label{lemma:suff:dec}

(\bfit{Sufficient Decrease Condition}) $F(\x^{t+1})  - F(\x^{t}) \leq  {- \frac{\theta}{2 g(\x^{t+1})} \|\x^{t+1} - \x^t\|_2^2 }$.

\end{lemma}
\noi \textbf{Remarks}. \textit{\textbf{(i)}} Both \textit{\textbf{FCD}} and \textit{\textbf{PCD}} have the same sufficient decrease condition. \textit{\textbf{(ii)}} The proximal parameter $\theta$ is critically important to guarantee global convergence and convergence rate of our algorithm.

\begin{lemma} \label{lemma:omega:sandwich}
(\bfit{Property of \textit{\textbf{FCD}}}) The value of the parameter $\alpha^t$ defined in (\ref{eq:opt:111}) is sandwiched as $F(\x^{t+1}) \leq \alpha^t \leq F(\x^{t+1}) + \sigma (F(\x^t) - F(\x^{t+1})) \leq  \sigma F(\x^0)$ with $\sigma \triangleq \frac{\max(\c)+\theta}{\theta}$.

\end{lemma}

We assume that the coordinate $i^t$ in each iteration is selected randomly and randomly. Our algorithm generates a random output $\x^t$ with $t=0,1,...$, which depends on the observed realization of the random variable: $\xi^{t-1} \triangleq \{i^0,~i^1,...,i^{t-1}\}$. We use $\E[\cdot]$ to denote the expectation of a random variable.

\subsection{Convex-Convex FMPs}\label{sebsect:convex}

This subsection presents some theoretical analysis for Algorithm \ref{algo:main} when the denominator $g(\cdot)$ is convex but not necessarily differentiable.

We first present the following useful definition.

\begin{definition} \label{ass:1}
(\textbf{Globally or Locally $\rho$-Bounded Non-Convexity}) \bfit{(i)} A function $\tilde{g}(\x)=-g(\x)$ is \textit{globally $\rho$-bounded non-convex} if: $\forall \x,\y,~\tilde{g}(\x) \leq \tilde{g}(\y)  + \la \partial \tilde{g}(\x),~\x - \y\ra + \frac{\rho}{2}\|\x - \y\|_2^2$ with $\rho<+\infty$. \bfit{(ii)} $\tilde{g}(\x)$ is \textit{locally $\rho$-bounded non-convex} if $\x$ is defined as some point $\check{\x}$ with $\x \triangleq \check{\x}$.

\end{definition}


\noi \textbf{Remarks}. \bfit{(i)} The definition of globally bounded non-convexity is also known as \textit{weakly-convex}, \textit{semi-convex}, or \textit{approximate convex} in the literature (cf. \cite{allen2018natasha,bohm2021variable,li2021weakly}). \bfit{(ii)} By this definition, $(\tilde{g}(\x)+\tfrac{\rho}{2}\|x\|_2^2)$ is convex. \bfit{(iii)} Smoothness is not required as convex functions are not necessarily smooth. \bfit{(iv)} It is not hard to verify that $\tilde{g}(\x) = - \|\G\x\|_4^2$ in (\ref{eq:app:eig:2}) is concave and globally bounded non-convex, while $ \tilde{g}(\x)= - \gamma \sum_{j=1}^k |\x_{[j]}|$ as in (\ref{eq:app:sparse}) is concave and locally bounded non-convex. See Section \ref{ass:globally:locally} in the \bfit{Appendix}.


%

\subsubsection{Optimality Analysis}


We now present two kinds of stationary solution which are novel in this paper.


\begin{definition}
(Fractional Coordinate-Wise Point, or \textit{FCW}-Point for short) Given a constant $\theta\geq 0$. Define $\mathcal{K}_i(\x,\eta)  \triangleq  \frac{\J_i(\x,\eta) }{g(\x + \eta\ei)}$. A solution $\ddot{\x}$ is called a \textit{FCW}-point if: $\mathcal{K}_i(\ddot{\x},0) = \min_{\etas_i}\mathcal{K}_i(\ddot{\x},\etas_i),~\forall i = 1,...,n$.

\end{definition}

\begin{definition}
(Parametric Coordinate-Wise Point, or \textit{PCW}-Point for short) Given a constant $\theta\geq 0$. Define $\mathcal{M}_i(\x,\eta)  \triangleq \J_i(\x,\eta)   - F(\x)   g(\x + \eta \ei)$. A solution $\dot{\x}$ is called a \textit{PCW}-point if: $\mathcal{M}_i(\dot{\x},0) = \min_{\etas_i}\mathcal{M}_i(\dot{\x},\etas_i),~\forall i = 1,...,n$.

\end{definition}

\noi \textbf{Remarks}. Both the \textit{FCW}-point and the \textit{PCW}-point use another \textit{non-convex} problem to characterize their stationary, and they state that if we minimize the majorization/surrogate function $\mathcal{K}_i(\ddot{\x},\eta)$ (or $\mathcal{M}_i(\dot{\x},\eta)$), we can not improve the objective value for $\mathcal{K}_i(\ddot{\x},\eta)$ (or $\mathcal{M}_i(\dot{\x},\eta)$) for all $i$.

\begin{lemma} \label{label:property}
For any \textit{FCW}-point $\ddot{\x}$ and any \textit{PCW}-point $\dot{\x}$, assume that $\tilde{g}(\x) = -g(\x)$ is locally $\rho$-bounded non-convex at the point $\ddot{\x}$ (or $\dot{\x}$) with $\rho<+\infty$. We define $\mathcal{C}(\x,\etas) \triangleq  \frac{1}{2}\| \etas \|^2_{\c+\theta} +  \frac{\rho}{2}\|\etas\|_2^2  F(\x)$. We have: $\textit{\textbf{(i)}}~\forall \etas,~F(\ddot{\x})  - F(\ddot{\x}+\etas)  \leq  \frac{\mathcal{C}(\ddot{\x},\etas)}{g(\ddot{\x}+\etas)}, ~\textit{\textbf{(ii)}}~\forall \etas,~F(\dot{\x})  - F(\dot{\x}+\etas)  \leq  \frac{\mathcal{C}(\dot{\x},\etas)}{g(\dot{\x}+\etas)}$.

\end{lemma}
\noi \textbf{Remarks}. The lemma above essentially implies that the optimality of \textit{FCW}-point coincides with that of \textit{PCW}-point; i.e., any \textit{FCW}-point must be a \textit{PCW}-point, and vice versa. 

We use $\check{\x}$, $\grave{\x}$, $\dot{\x}$, $\ddot{\x}$, and $\bar{\x}$ to denote a \textit{C}-point, a \textit{D}-point, a \textit{FCW}-point, a \textit{PCW}-point, and an optimal point, respectively. The following theorem establishes their relations.

\begin{theorem} \label{the:optimality}

\textbf{(Optimality Hierarchy between the Optimality Conditions).} Based on the the assumption made in Lemma \ref{label:property}. The following relations hold: $\{\bar{\x}\} \overset{\textbf{(a)}}{ \subseteq}\{\ddot{\x}\}  \overset{\textbf{(b)}}{\Leftrightarrow } \{\dot{\x}\} \overset{\textbf{(c)}} { \subseteq }  \{\grave{\x}\} \overset{\textbf{(d)}} { \subseteq }  \{\check{\x}\}$.

\end{theorem}

\noi \textbf{Remarks}. The optimality condition of \textit{FCW}-point or \textit{PCW}-point is stronger than that of \textit{C}-point \cite{LIHarmonic,LISIOPT2022,Radu2017} and \textit{D}-point \cite{PangRA17} when $(-g(\cdot))$ is $\rho$-bounded non-convex. We use the following one-dimensional example to clarify this point: $\min_{x} f(\x)  \triangleq \frac{(x+2)^2}{|3x+2|+1}$. This problem contains three \textit{C}-points $\{-2,-\tfrac{2}{3},0\}$ and two \textit{D}-points $\{-2,0\}$. $x=-2$ is the unique \textit{FCW}-point since it is the unique global optimal solution for this one-dimensional problem. After some preliminary calculations, one can verify that $x=-2$ is also the unique \textit{PCW}-point. See Section \ref{app:sebsect:hierarchy} in the \bfit{Appendix}. 



\subsubsection{Convergence Analysis}

We first define the approximate \textit{FCW}-Point and approximate \textit{PCW}-Point.
\begin{definition}
(Approximate \textit{FCW}-Point and Approximate \textit{PCW}-Point) Given any constant $ \epsilon>0$. Define $\J_i(\x,\eta)  \triangleq f(\x) +  \nabla_i f(\x) \eta + \frac{\c_i+\theta}{2}\eta^2 + h(\x+\eta\ei)$ with $\theta\geq 0$, $\mathcal{K}_i(\x,\eta)  \triangleq  \frac{\J_i(\x,\eta) }{g(\x + \eta\ei)}$, $\mathcal{M}_i(\x,\eta)  \triangleq  \frac{\J_i(\x,\eta) }{g(\x + \eta\ei)}$.

~\textit{\textbf{(i)}} A solution $\ddot{\x}$ is called an $\epsilon$-approximate \textit{FCW}-point if: $\frac{1}{n}\sum_{i=1}^n \text{dist}(0,\arg \min_{\eta} \mathcal{K}_{i}(\ddot{\x},\eta)   )^2 \leq \epsilon$.

\textit{\textbf{(ii)}} A solution $\dot{\x}$ is called an $\epsilon$-approximate \textit{PCW}-point if: $\frac{1}{n}\sum_{i=1}^n \text{dist}(0,\arg \min_{\eta} \mathcal{M}_{i}(\dot{\x},\eta)   )^2 \leq \epsilon$.
\end{definition}

%
%

We now prove the global convergence of Algorithm \ref{algo:main} for convex-convex FMPs.

%

\begin{proposition} \label{prop:convex:conv}
(Global Convergence) Assume that $\x^t$ is bounded for all $t$\footnote{This condition always holds if we impose bound constraints on $\x$ that $\|\x\|_{\infty}\leq \vartheta$, refer to Problem (\ref{eq:app:sparse}).}, any clustering point of the sequence is almost surely a \textit{FCW}-point (or a \textit{PCW}-point) of Problem (\ref{eq:main}). Furthermore, Algorithm \ref{algo:main} finds an $\epsilon$-approximate \textit{FCW}-point (or \textit{PCW}-point) of Problem (\ref{eq:main}) in at most $T+1$ iterations in the sense of expectation, where $T\leq \lceil\tfrac{2 n \bar{g} F(\x^0) }{\theta \epsilon }\rceil = \mathcal{O}(\epsilon^{-1})$.

\end{proposition}

To achieve stronger convergence result for Algorithm \ref{algo:main}, we make the following additional assumption.

\begin{assumption}\label{ass:local:tsneg:luo:bound}
(Luo-Tseng Error Bound \cite{luo1993error,tseng2009coordinate}) We define a residual function as $\mathcal{R}(\x) \triangleq \frac{1}{n} \sum_{i=1}^n |\P_i(\x)|$, where $\P_i(\x)$ is defined in (\ref{eq:subprob:nonconvex}) (or (\ref{eq:subprob:nonconvex2})). For any $\zeta\geq \min_{\x}F(\x)$, there exist scalars $\delta>0$ and $\varrho>0$ such that:
\beq \label{eq:luotseng}
 \text{dist}(\x,\mathcal{X}) \leq \delta \cdot \R(\x) ,~\text{whenever} ~F(\x)\leq \zeta, ~\R(\x) \leq \varrho.
\eeq
\noi Here, $dist(\x,\mathcal{X}) = \inf_{\z\in \mathcal{X}}\|\z-\x\|$, $\mathcal{X}$ is the set of the \textit{FCW}-point (or the \textit{PCW}-point).

\end{assumption}

Luo-Tseng error bound has long been a significant topic in all aspects of mathematical optimization. Many optimization problems have been shown to possess the Luo-Tseng error bound property \cite{YueZS19,Dong2021}. Assumption \ref{ass:local:tsneg:luo:bound} is similar to the classical local proximal error bound assumption in the literature. We note that if $\x^t$ is not the \textit{FCW}-point (or the \textit{PCW}-point), we have $\R(\x^t)>0$. By the boundedness of $\x^t$ and $\dot{\x}$ (or $\ddot{\x}$), there exists a sufficiently large constant $\delta$ such that $\text{dist}(\x^t,\mathcal{X}) \leq \delta \cdot \R(\x^t)$. Thus, Assumption \ref{ass:local:tsneg:luo:bound} is reasonable.

We now establish the convergence rate for Algorithm \ref{algo:main}. We have the following two theorems.

\begin{theorem}\label{the:rate:FCDC}
{(Convergence Rate of \textit{\textbf{FCD}})}. For any \textit{FCW}-point $\ddot{\x}$, we define $\ddot{q}^t \triangleq F(\x^t) - F(\ddot{\x}),~r^t  \triangleq \frac{1}{2}\|\x^t - \ddot{\x}\|_{\bar{\c}}^2,~\bar{\c} \triangleq \c + \theta$.  Assume that $\tilde{g}(\x)=-g(\x)$ is globally $\rho$-bounded non-convex, and $F(\cdot)$ satisfies Assumption \ref{ass:local:tsneg:luo:bound}. We define: $\varpi \triangleq \frac{\max(\bar{\c}) }{\min(\bar{\c})} \cdot\frac{\rho}{\theta} \cdot F(\x^0)$. We have the following inequality: $(1- \varpi) \E_{i^{t}}[r^{t+1}]  + \frac{g(\bar{\x})}{n}\ddot{q}^{t+1} \leq (1-\varpi) r^t + \frac{\varpi}{n} r^t$. When the proximal parameter $\theta$ is sufficiently large such that $\varpi\leq 1$, we obtain: $\ddot{q}^{t+1} \leq (\frac{\kappa_1}{\kappa_1 + \kappa_0})^{t+1} \ddot{q}^0$, where $\kappa_0 \triangleq\frac{g(\bar{\x})}{\bar{g}}$ and $\kappa_1\triangleq(n + 1) \max(\bar{\c})\delta^2/\theta$.

\end{theorem}


\begin{theorem}\label{the:rate:DCDC}
{(Convergence Rate of \textit{\textbf{PCD}})}. For any \textit{PCW}-point $\dot{\x}$, we define $\dot{q}^t \triangleq F(\x^t) - F(\dot{\x}),~r^t  \triangleq \frac{1}{2}\|\x^t - \dot{\x}\|_{\bar{\c}}^2,~\bar{\c} \triangleq \c + \theta$. Assume that $\tilde{g}(\x)=-g(\x)$ is globally $\rho$-bounded non-convex, and $F(\cdot)$ satisfies Assumption \ref{ass:local:tsneg:luo:bound}. We define: $\varpi \triangleq \frac{\rho}{\min(\bar{\c})}F(\x^0) $. We have the following inequality: $\E_{i^{t}}[(1-\varpi) r^{t+1} ] + \frac{\bar{g}}{n} \dot{q}^{t+1} \leq (1-\varpi) r^t + \frac{\varpi}{n}r^t - \frac{g(\bar{\x})}{n}\dot{q}^t + \frac{\bar{g}}{n} \dot{q}^t$. When the proximal parameter $\theta$ is sufficiently large such that $\varpi \leq1$, we obtain: $\dot{q}^{t+1} \leq   (\frac{\kappa_1+1 - \kappa_0}{\kappa_1 + 1})^{t+1}\dot{q}^0$, where $\kappa_0\triangleq\frac{g(\bar{\x})}{\bar{g}}$ and $\kappa_1 \triangleq(n + 1) \max(\bar{\c})\delta^2/\theta$.

\end{theorem}

\noi \textbf{Remarks}. \textit{\textbf{(i)}} Algorithm \ref{algo:main} converges to the \textit{FCW}-point (or the \textit{PCW}-point) with a Q-linear convergence rate. \textit{\textbf{(ii)}} We compare the convergence rate of \textit{\textbf{FCD}} and \textit{\textbf{PCD}} which depend on $\kappa_0$ and $\kappa_1$: $(\frac{\kappa_1+1 - \kappa_0}{\kappa_1 + 1})-(\frac{\kappa_1}{\kappa_1 + \kappa_0}) = \frac{ 1  }{(\kappa_1 + \kappa_0)(\kappa_1 + 1)} [ \kappa_1(\kappa_1 + \kappa_0) +  (\kappa_1 + \kappa_0) -  \kappa_0(\kappa_1 + \kappa_0) - \kappa_1 (\kappa_1 + 1) ] = \frac{  \kappa_0 (1-\kappa_0)}{(\kappa_1 + \kappa_0)(\kappa_1 + 1)} \geq 0$, where the inequality holds due to $\kappa_0\triangleq\frac{g(\bar{\x})}{\bar{g}} \in (0,1]$. Thus, the convergence rate of $\textit{\textbf{FCD}}$ is better than that of \textit{\textbf{PCD}} in theory.





\subsection{Convex-Concave FMPs}\label{sebsect:concave}

This subsection provides some theoretical analysis for Algorithm \ref{algo:main} when the donominator $g(\cdot)$ is concave and differentiable \footnote{\label{foot1}Since any coordinate-wise stationary point is not necessarily the first-order stationary point for non-separable and non-differentiable convex functions \cite{tseng2009coordinate} as in (\ref{eq:subprob:nonconvex2}), we further assume that the convex function $(-g(\cdot))$ is differentiable.}. Convex-concave FMPs as in Problem \ref{algo:main} can be converted to an equivalent convex program using the Charnes-Cooper transformation \cite{hadjisavvas2006handbook}: $\min_{t,\y} t f(\y/t) + t h(\y/t),~s.t.~tg(\y/t) \geq 1$. However, our CD methods are able to directly solve Problem (\ref{eq:main}) with exploiting its specific structure.

$\bullet$ \textbf{Optimality Analysis.} We provide some theoretical insights into convex-concave FMP.

\begin{proposition} \label{proposition:global}
(i) $F(\cdot)$ is quasiconvex that: $F( \alpha \x + (1-\alpha) \y) \leq \max(F(\x),F(\y)), \forall \alpha \in [0,1], \x,\y.$ (ii) Any critical point of Problem (\ref{eq:main}) is a global minimum.

\end{proposition}
\noi \textbf{Remarks}. \bfit{(i)} Using an inequality $\forall a\geq 0,b\geq 0,c>0,d>0,~\frac{a + b}{c+d} \leq \max(\frac{a}{c},\frac{b}{d})$, we prove the quasiconvexity of $F(\cdot)$. \bfit{(ii)} It is shown that any local minimum of a \textit{strictly quasiconvex problem} is also a global minimum (cf. section 3.5.5 in \citet{bazaraa2013nonlinear}). General quasiconvex problems do not enjoy this property while we prove that convex-concave FMPs do.






$\bullet$ \textbf{Convergence Analysis.} We now establish the \textit{convergence rate} of Algorithm \ref{algo:main}. 

\begin{theorem}\label{the:rate:FCD:PCD:concave}
{(Convergence Rate)}. For any global optimal solution $\bar{\x}$-point of Problem (\ref{eq:main}), we define $q^t \triangleq F(\x^t) - F(\bar{\x}),~r^t  \triangleq \frac{1}{2}\|\x^t - \bar{\x}\|_{\bar{\c}}^2,~\bar{\c} \triangleq \c + \theta$. \bfit{(i)} For \textit{\textbf{FCD}}, we have: $\E_{\xi^{t-1}}[q^t] \leq  \frac{n(\bar{g} \sigma  q^0  + r^0)}{  g(\bar{\x})t}$, where $\sigma$ is defined in Lemma \ref{lemma:omega:sandwich}. \bfit{(ii)} For \textit{\textbf{PCD}}, we have: $\E_{\xi^{t-1}}[q^t] \leq  \frac{n(\bar{g} q^0 + r^0)}{g(\bar{\x}) (t+1) }$.
\end{theorem}

\noi \textbf{Remarks}. Algorithm \ref{algo:main} converges to the global optimal solutions with a sublinear convergence rate.

\section{Implementations and Experiments} \label{sect:exp}

We first describe the implementations of Algorithm \ref{algo:main} for solving the sparse recovery problem and the $\ell_p$ norm eigenvalue problem, and then provide numerical comparisons against state-of-the-art methods on some real-world data. Since the two CD methods achieve the same optimality condition (refer to Theorem \ref{the:optimality}) and we pay more attention to better optimality/accuracy, we only implement one of them for comparisons.

\subsection{Implementations for Sparse Recovery}

We observe that Problem (\ref{eq:app:sparse}) is a special case of Problem (\ref{eq:main}) with $f(\x) = \frac{1}{2}\| \G \x  - \y \|_2^2$, $h(\x) = \gamma\|\x\|_1 + I_{\Delta}(\x)$, and $g(\x) = \gamma \sum_{j=1}^k |\x|_{[j]}$, where $I_{\Delta}(\x)={\tiny \{
                               \begin{array}{ll}
                                 0,  & \hbox{ $\x \in \Delta$;} \\
                                 +\infty, & \hbox{else.}
                               \end{array}
                             \}}$, $\Delta\triangleq \{\x|\|\x\|_{\infty}\leq \vartheta\}$. The gradient of $f(\x)$ can be computed as: $\nabla f(\x)=\G^T (\G\x - \y)\triangleq \g$. $\nabla f(\x)$ is $L$-Lipschitz continuous with $L=\|\G\|_2^2$ and coordinate-wise Lipschitz with $\c_i = (\G\G^T)_{ii},~\forall i$. The subgradient of $g(\x)$ can be computed as $(\partial g(\x))_i ={\tiny \{
                               \begin{array}{ll}
                                 \text{sign}(\x_i),  & \hbox{ $i \in \Delta_k(\x)\text{~and~} \x_i \neq 0$;} \\
                                 $[-1,1]$, & \hbox{else.}
                               \end{array}
                             \}}$, where $\Delta_k(\x)$ is the index of the largest (in magnitude) $k$ elements of $\x$. According to Algorithm \ref{algo:main}, the update for \textit{\textbf{PCD}} reduces to solving the following one-dimensional problem: $\textstyle \min_{c_1\leq \eta\leq c_2}  \frac{a}{2} \eta^2    + \eta b + \gamma \|\x + \eta \ei\|_1 - \tau     \sum_{j=1}^k |\x + \eta \ei|_{[j]}$, where $a=\c_{i}+\theta$,~$b=\g_i$, $\tau=\gamma F(\x^t)$, $c_1=-\vartheta-\x_i$, $c_2=\vartheta-\x_i$. We choose $\textit{\textbf{PCD}}$ for comparisons since its subproblem is easier to solve. 

$\bullet$ \textbf{A Breakpoint Searching Procedure for} \textit{\textbf{PCD}}. At first, we drop the bound constraint $c_1\leq \eta\leq c_2$. Since the variable $\eta$ only affects the value of $\x_{i}$, we consider two cases for $\x_{i}+\eta$. \textit{\textbf{(i)}} $\x_{i}+\eta$ belongs to the top-$k$ subset. Problem (\ref{eq:subprob:nonconvex2}) reduces to $\min_{\eta} \frac{a}{2} \eta^2 + \eta b +  \gamma |\x_{i} + \eta | - \tau | \x_{i} + \eta |$. We consider three cases for the non-smooth term $|\x_{i} + \eta |$, leading to three breakpoints: $\{-\x_{i}, (\tau      - \gamma - b)/a,(\gamma - \tau - b)/a\}$. \textit{\textbf{(ii)}} $\x_i+\eta$ does not belong to the top-$k$ subset. Problem (\ref{eq:subprob:nonconvex2}) reduces to $\min_{\eta} \eta b + \frac{a}{2} \eta^2 + \gamma |\x_{i} + \eta |$. Again, we consider three cases for the term $|\x_{i} + \eta |$, resulting in three breakpoints: $\{-\x_{i^t},(-\gamma-b)/a,~(\gamma-b)/a\}$. Therefore, Problem (\ref{eq:subprob:nonconvex2}) contains $5$ different breakpoints $\Theta'=\{-\x_{i}, (\tau      - \gamma - b)/a,(\gamma - \tau - b)/a,(-\gamma-b)/a,~(\gamma-b)/a\}$ without the bound constraint. At last, taking the bound constraint into consideration, we conclude that Problem (\ref{eq:subprob:nonconvex2}) contains 7 breakpoints $\Theta = \{c_1,c_2, \min(c_2,\max(c_1,\Theta'))\}$.

Once we have identified all the possible breakpoints / critical points $\Theta$ for the one-dimensional subproblem, we pick the solution that leads to the lowest value as the optimal solution.

$\bullet$ \textbf{Compared Methods}. We compare \bfit{PCD} against the following three methods. For ease of discussion, we only consider $\vartheta=+\infty$ in the sequel. \bfit{(i)} \bfit{DPA} iteratively generates a sequence $\{\x^t\}$ as: $\x^{t+1} = \arg \min_{\x}f(\x) + \gamma\|\x\|_1 - \lambda^t \la\x -\x^t,\partial g(\x^t) \ra$, which is solved by an accelerated proximal gradient method \cite{beck2009fast,nesterov2003introductory}. \bfit{(ii)} \bfit{PGSA} generates the new iterate by: $\x^{t+1} = \arg \min_{\x} \frac{L}{2}\|\x -\x^t\|_2^2+ \la \x-\x^t, \nabla f(\x^t)\ra + f(\x^t) + \gamma\|\x\|_1 - F(\x^t) \la\x -\x^t,\partial g(\x^t) \ra$, which reduces a soft-thresholding operator. \bfit{(iii)} Quadratic Transform Parametric Algorithm (QTPA) solves the following problem: $\min_{\x}~ f(\x) + h(\x) -  2 (\beta^t)^{-1} \sqrt{ g(\x) }$ with $\beta^t$ is renewed as: $\beta^t= {\sqrt{g(\x^t)}}/ { [f(\x^t) + h(\x^t)]}$. We consider a proximal gradient-subgradient algorithm \cite{LISIOPT2022,boct2021extrapolated} to minimize the objective function over $\x$, leading to the following update: $
\x^{t+1} = \arg \min_{\x}\frac{L}{2}\|\x -\x^t\|_2^2+  \la \x-\x^t, \nabla f(\x^t)\ra     +f(\x^t) +  \gamma\|\x\|_1 -  2 (\beta^t)^{-1} \la\x -\x^t,\partial \breve{g}(\x^t) \ra$, where $\partial \breve{g}(\x^t)$ is the subgradient of $\sqrt{ g(\x) }$ which can be computed as $\partial\breve{g}(\x) = \frac{1}{2}g(\x)^{-1/2}\cdot \partial {g}(\x)$.

\subsection{Implementations for $\ell_p$ Norm Eigenvalue Problem}

We observe that Problem (\ref{eq:app:eig:2}) is a special case of Problem (\ref{eq:main}) with $f(\x) = \x^T\mathbf{Q}\x +\gamma_3$, $h(\x) =  0 $, and $g(\x) = \|\G\x\|^2_p + \gamma_4$. We consider the classical ICA problem and choose $p=4$, $\mathbf{Q}=\mathbf{I}$. We simply set $\gamma_3=\gamma_4=0$. Therefore, we have $F(\x) = \frac{\|\x\|_2^2}{\|\G\x\|_4^2}$. The gradient of $f(\x)$ can be computed as $\nabla f(\x)=2\x $. $\nabla f(\x)$ is coordinate-wise Lipschitz continuous with $\c_i=2,~\forall i$. The gradient of $\g(\x)$ can be computed as $\nabla g(\x) = 2 g(\x) \cdot \sum_{i=1}^m  ( (\G_i\x)^3 \G_i^T )$ with $\G_i \in \mathbb{R}^{1\times n}$ being the $i$-the row of $\G$. According to Algorithm \ref{algo:main}, the update for \textit{\textbf{FCD}} reduces to solving the following one-dimensional problem: $\min_{\eta}\frac{ \|\x^t\|_2^2 + 2 \x_i \eta + \frac{2+\theta}{2} \eta^2 }{ \sqrt{\|\G(\x^t + \eta \eit)\|^4_4}}$. We choose \textit{\textbf{FCD}} for comparisons since its subproblem is easier to solve.

$\bullet$ \textbf{A Breakpoint Searching Procedure for} \textit{\textbf{FCD}}. We note that one-dimensional problem boils down to the following problem: $\min_\eta p(\eta) \triangleq \frac{ a_2 \eta^2 + a_1 \eta + a_0}{ \sqrt{ b_4 \eta^4 + b_3 \eta^3+ b_2 \eta^2 + b_1 \eta + b_0 }}$ with suitable parameters $a_2,a_1,a_0$ and $b_4,b_3,b_2,b_1,b_0$. Setting the gradient of $p(\cdot)$ to zero yields: $2 a_2 \eta + a_1 = p(\eta)  \tfrac{1}{2}(b_4 \eta^4 + b_3 \eta^3 + b_2 \eta^2 + b_1 \eta + b_0)^{-\frac{1}{2}} \cdot (4 b_4 \eta^3 + 3 b_3 \eta^2 + 2 b_2 \eta + b_1)$. After some preliminary calculations, the equation above is equivalent to the following quartic equation: $c_4 \eta^4 + c_3 \eta^3 + c_2 \eta^2 +  c_1 \eta + c_0  = 0$ with suitable parameters $c_4,c_3,c_2,c_1,c_0$. It can be solved analytically by Lodovico Ferrari's method (\url{https://en.wikipedia.org/wiki/Quartic_equation}).

$\bullet$ \textbf{Compared Methods}. We compare \textit{\textbf{PCD}} against the following two methods. \textit{\textbf{(i)}} The power Method \cite{hyvarinen2000independent} solves the original problem $\max_{\v} \|\G\v\|^4_4,~s.t.~\|\v\| = 1$ using the following update: $\x^{t+1} =   \frac{\partial  {g}(\x^t)}{ \|\partial  {g}(\x^t)\|}$. \textit{\textbf{(i)}} \textbf{PGSA} generates the new iterate by: $\x^{t+1} = \arg \min_{\x}~\|\x\|_2^2 - F(\x^t) \la \x - \x^t,~\partial g(\x^t) \ra = \tfrac{1}{2} F(\x^t) \partial g(\x^t)$. Interestingly, we find that the solution of PGSA and that of the power method only differ by a scale factor. Since the objective function $F(\x)\triangleq \frac{\|\x\|_2^2}{\|\G\x\|^4_2}$ is scale invariance, these two solutions lead to the same objective value for all iterations.

\begin{table}[!t]
\centering
\scalebox{0.63}{\begin{tabular}{|c|c|c|c|c|c|c|c|c|}
\hline
           & \textit{\textbf{DPA}}  & \textit{\textbf{PGSA}} & \textit{\textbf{QTPA}} & \textit{\textbf{PCD}}   \\
\hline
 e2006-1000-1024 & 1.874 $\pm$ 0.315 & 1.929 $\pm$ 0.278 & 1.923 $\pm$ 0.279 & \cone{1.530 $\pm$ 0.184}  \\
 e2006-1000-2048 & 1.640 $\pm$ 0.118 & 1.663 $\pm$ 0.172 & 1.660 $\pm$ 0.177 & \cone{1.312 $\pm$ 0.061}  \\
 e2006-1024-1000 & 2.610 $\pm$ 0.796 & 2.362 $\pm$ 0.533 & 2.362 $\pm$ 0.530 & \cone{1.882 $\pm$ 0.418}  \\
 e2006-2048-1000 & 5.623 $\pm$ 4.005 & 6.576 $\pm$ 4.966 & 6.593 $\pm$ 4.989 & \cone{3.068 $\pm$ 1.282}  \\
 news20-1000-1024 & 1.750 $\pm$ 0.247 & 1.403 $\pm$ 0.128 & 1.402 $\pm$ 0.130 & \cone{1.168 $\pm$ 0.023}  \\
 news20-1000-2048 & 2.043 $\pm$ 0.429 & 1.424 $\pm$ 0.181 & 1.426 $\pm$ 0.180 & \cone{1.207 $\pm$ 0.065}  \\
 news20-1024-1000 & 1.856 $\pm$ 0.353 & 1.488 $\pm$ 0.317 & 1.487 $\pm$ 0.318 & \cone{1.195 $\pm$ 0.045}  \\
 news20-2048-1000 & 4.997 $\pm$ 0.269 & 2.664 $\pm$ 0.604 & 2.559 $\pm$ 0.745 & \cone{1.394 $\pm$ 0.115}  \\
 sector-1000-1024 & 1.864 $\pm$ 0.162 & 1.337 $\pm$ 0.105 & 1.337 $\pm$ 0.104 & \cone{1.160 $\pm$ 0.016}  \\
 sector-1000-2048 & 1.780 $\pm$ 0.040 & 1.293 $\pm$ 0.033 & 1.293 $\pm$ 0.026 & \cone{1.148 $\pm$ 0.010}  \\
 sector-1024-1000 & 2.039 $\pm$ 0.016 & 1.485 $\pm$ 0.194 & 1.486 $\pm$ 0.195 & \cone{1.193 $\pm$ 0.015}  \\
 sector-2048-1000 & 5.041 $\pm$ 1.714 & 2.477 $\pm$ 1.048 & 2.475 $\pm$ 1.046 & \cone{1.409 $\pm$ 0.108}  \\
 TDT2-1000-1024 & 1.778 $\pm$ 0.303 & 1.646 $\pm$ 0.035 & 1.644 $\pm$ 0.032 & \cone{1.215 $\pm$ 0.047}  \\
 TDT2-1000-2048 & 1.710 $\pm$ 0.045 & 1.398 $\pm$ 0.029 & 1.398 $\pm$ 0.028 & \cone{1.127 $\pm$ 0.016}  \\
 TDT2-1024-1000 & 1.984 $\pm$ 0.284 & 1.555 $\pm$ 0.058 & 1.552 $\pm$ 0.050 & \cone{1.206 $\pm$ 0.067}  \\
 TDT2-2048-1000 & 4.696 $\pm$ 1.980 & 3.846 $\pm$ 0.901 & 3.789 $\pm$ 0.800 & \cone{1.338 $\pm$ 0.038}  \\
\hline
\end{tabular}}
\caption{Comparisons of objective values for solving the spare recovery problem.}
\label{tab:acc:sparse:recovery}
\end{table}

\vspace{-12pt}

\begin{figure}
\centering     
\subfigure[\scriptsize e2006-1000-2048]{\label{fig:a}\includegraphics[width=0.234\textwidth]{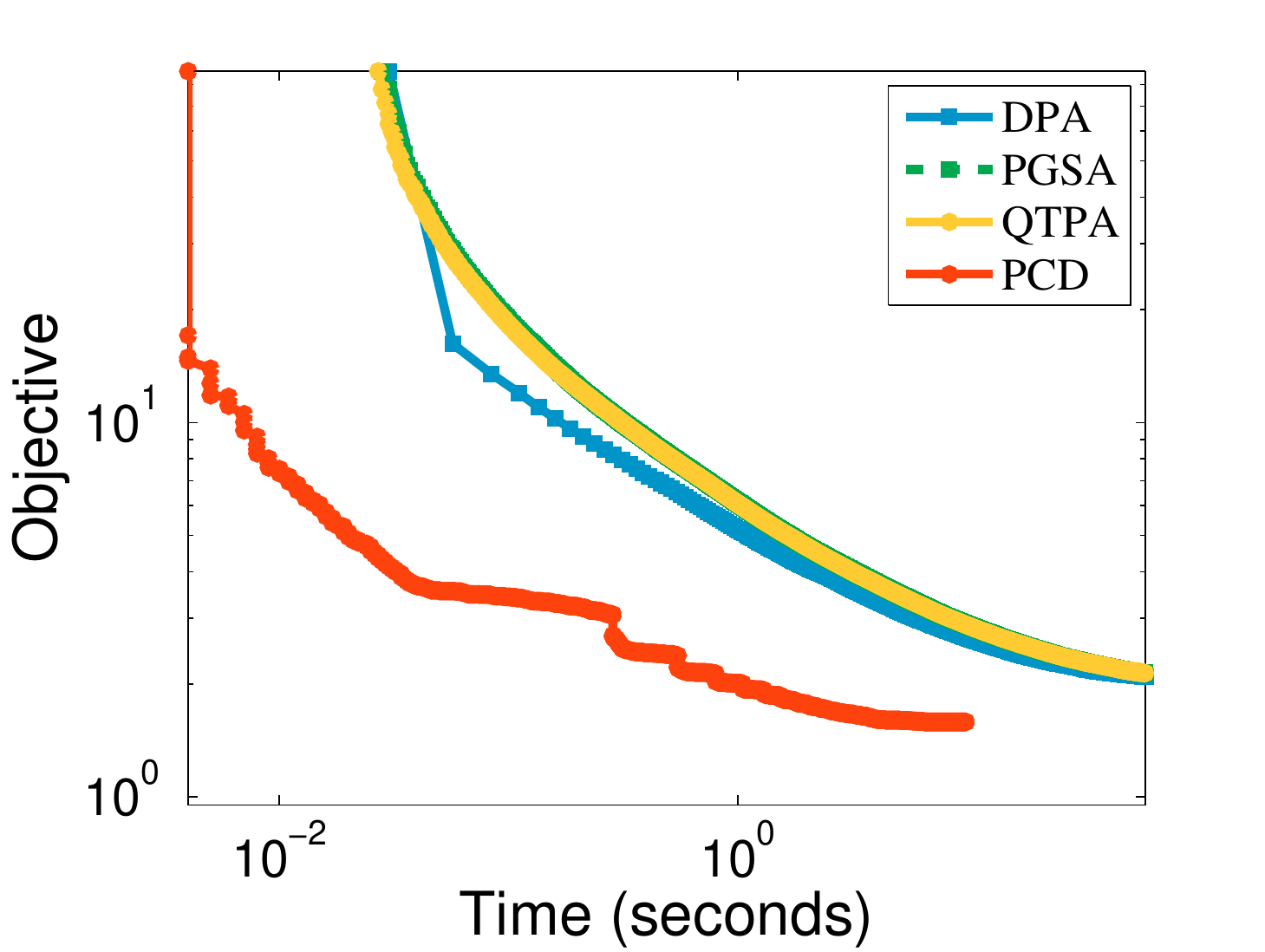}}
\subfigure[\scriptsize e2006-2048-1000]{\label{fig:b}\includegraphics[width=0.234\textwidth]{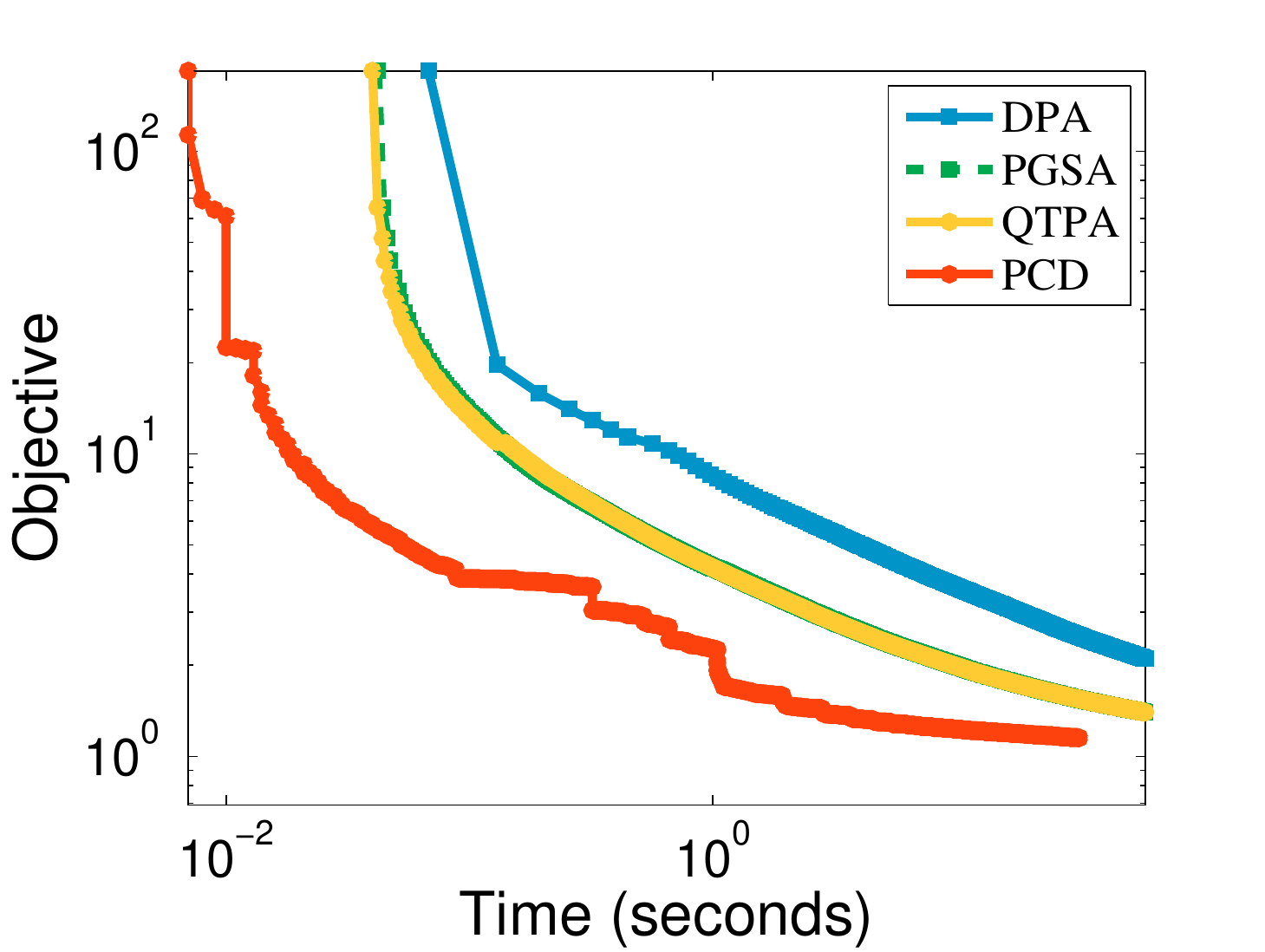}}
\caption{The convergence curve for solving the sparse recovery problem.} \label{figure:sparse}
\end{figure}

%
%
%

%
%

\subsection{Experiment Settings}

To generate the design/signal matrix $\G$, we consider four publicly available real-world data sets: `e2006tfidf', `news20', `sector', and `TDT2'. We randomly select a subset of examples from the original data sets (\url{http://www.cad.zju.edu.cn/home/dengcai/Data/TextData.html},
\url{https://www.csie.ntu.edu.tw/~cjlin/libsvm/}). The size of $\G\in \mathbb{R}^{m\times n}$ are are chosen from the following set $(m,n) \in \{ (1000,1024), (1000,2048), (1024, 1000), (2048, 1000)\}$. To generate the original $k$-sparse signal $\bar{\x}$ for the sparse recovery problem, we randomly select a support set $S$ of size $100$ and set $\bar{\x}_{\{1,...,n\} \setminus S}=\mathbf{0}$,~$\bar{\x}_{S} =\text{randn}(|S|,1)$. We generate the observation vector via $\y = \G\bar{\x} +  0.1 \|\G\bar{\x}\| \cdot \text{randn}(m,1)$. All methods are implemented in MATLAB on an Intel 2.6 GHz CPU with 64 GB RAM. We use the Matlab inbuilt function `roots' to solve the quartic equation. We define $\w_t  \triangleq [{F}(\x^t)-{F}(\x^{t+1})]/\max(1,{F}(\x^t))$, and let all algorithms run up to $T$ seconds and stop them at iteration $t$ if $\text{mean}([{\w}_{t-\text{min}(t,\upsilon)+1},\w_{t-min(t,\upsilon)+2},...,\z_t]) \leq \epsilon$. We use the default value $(\theta,\epsilon,\upsilon,T)=(10^{-6},10^{-10},500,100)$. All methods are executed 10 times and the average performance is reported. We only use the cyclic order rule to select the coordinate for Algorithm \ref{algo:main}. Some Matlab code can be found in the \bfit{supplemental material}.

\begin{table}[!t]
\centering
\scalebox{0.632}{\begin{tabular}{|c|c|c|c|c|c|c|c|c|}
\hline
            & \textit{\textbf{PGSA}} & \textit{\textbf{Power Method}} & \textit{\textbf{FCD}}    \\
\hline
e2006-1000-1024 & 12.254 $\pm$ 14.922 & 12.254 $\pm$ 14.922 & \cone{6.686 $\pm$ 4.956}  \\
e2006-1000-2048 & 16.896 $\pm$ 14.521 & 16.896 $\pm$ 14.521 & \cone{9.436 $\pm$ 6.359}  \\
e2006-1024-1000 & 5.923 $\pm$ 4.485 & 5.923 $\pm$ 4.485 & \cone{4.948 $\pm$ 2.631}  \\
e2006-2048-1000 & 16.846 $\pm$ 13.916 & 16.846 $\pm$ 13.916 & \cone{11.360 $\pm$ 8.225}  \\
news20-1000-1024 & 112.805 $\pm$ 58.995 & 112.805 $\pm$ 58.995 & \cone{78.183 $\pm$ 22.830 } \\
news20-1000-2048 & 125.440 $\pm$ 43.203 & 125.440 $\pm$ 43.203 & \cone{120.046 $\pm$ 41.353}  \\
news20-1024-1000 & 99.211 $\pm$ 35.338 & 99.211 $\pm$ 35.338 & \cone{80.244 $\pm$ 22.771}  \\
news20-2048-1000 & 138.909 $\pm$ 49.626 & 138.909 $\pm$ 49.626 & \cone{108.080 $\pm$ 37.811}  \\
sector-1000-1024 & 60.813 $\pm$ 24.018 & 60.813 $\pm$ 24.018 & \cone{50.551 $\pm$ 18.675}  \\
sector-1000-2048 & 139.459 $\pm$ 51.094 & 139.459 $\pm$ 51.094 & \cone{96.301 $\pm$ 42.115}  \\
sector-1024-1000 & 83.176 $\pm$ 38.697 & 83.176 $\pm$ 38.697 & \cone{48.559 $\pm$ 19.163}  \\
sector-2048-1000 & 104.654 $\pm$ 63.318 & 104.654 $\pm$ 63.318 & \cone{78.110 $\pm$ 28.532}  \\
TDT2-1000-1024 & 27.167 $\pm$ 12.705 & 27.167 $\pm$ 12.705 & \cone{22.308 $\pm$ 8.171}  \\
TDT2-1000-2048 & 27.480 $\pm$ 15.468 & 27.480 $\pm$ 15.468 & \cone{23.225 $\pm$ 12.614}  \\
TDT2-1024-1000 & 32.334 $\pm$ 18.178 & 32.334 $\pm$ 18.178 & \cone{21.143 $\pm$ 12.143}  \\
TDT2-2048-1000 & 44.659 $\pm$ 19.775 & 44.659 $\pm$ 19.775 & \cone{36.517 $\pm$ 12.689}  \\
\hline
\end{tabular}}
\caption{Comparisons of objective values for solving the $\ell_p$ Norm Eigenvalue Problem with $p=4$.} \label{tab:acc:ICA}
\end{table}

\vspace{-12pt}
\begin{figure}
\centering     
\subfigure[\scriptsize e2006-1000-2048]{\label{fig:a}\includegraphics[width=0.234\textwidth]{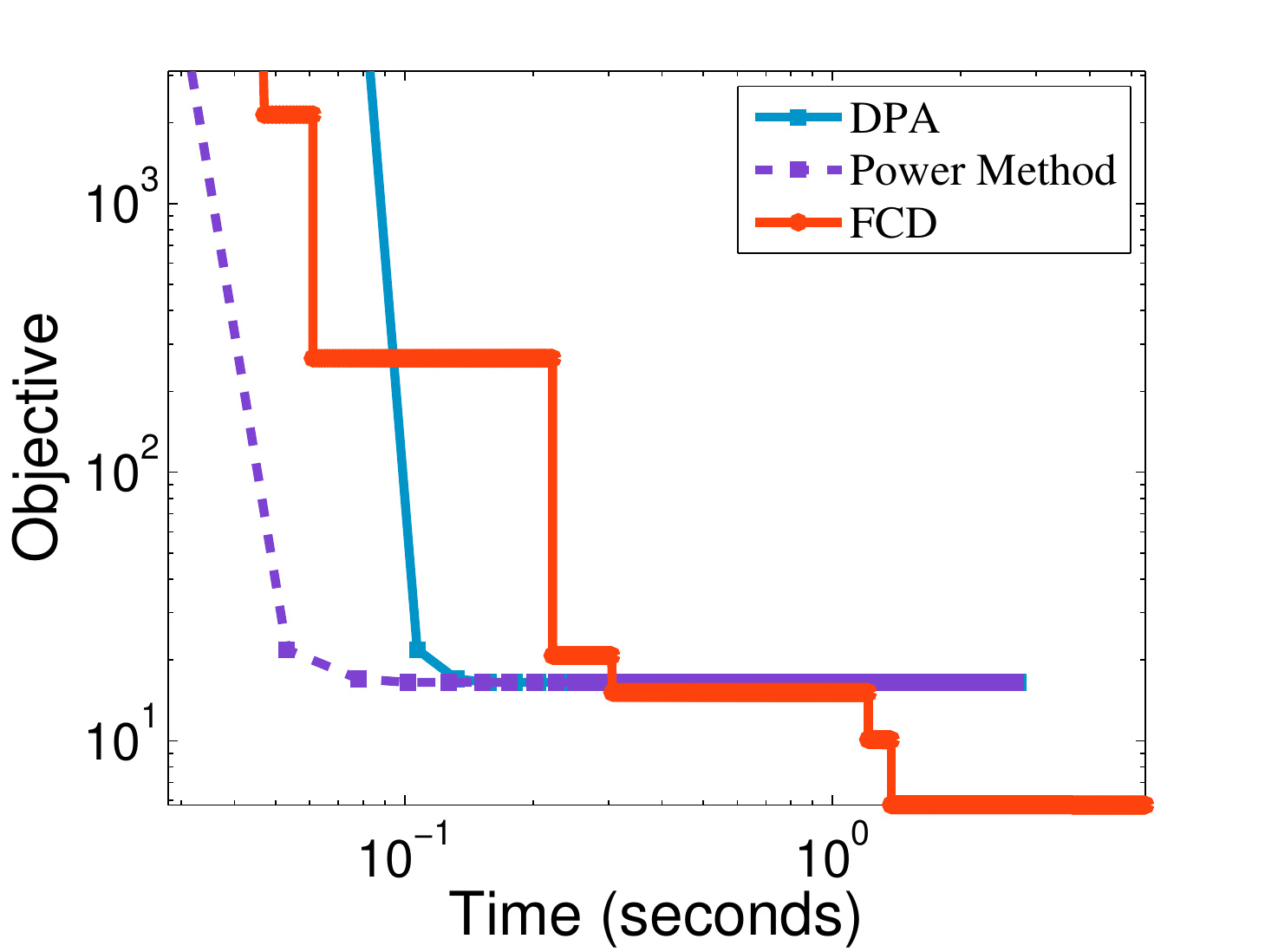}}
\subfigure[\scriptsize e2006-2048-1000]{\label{fig:b}\includegraphics[width=0.234\textwidth]{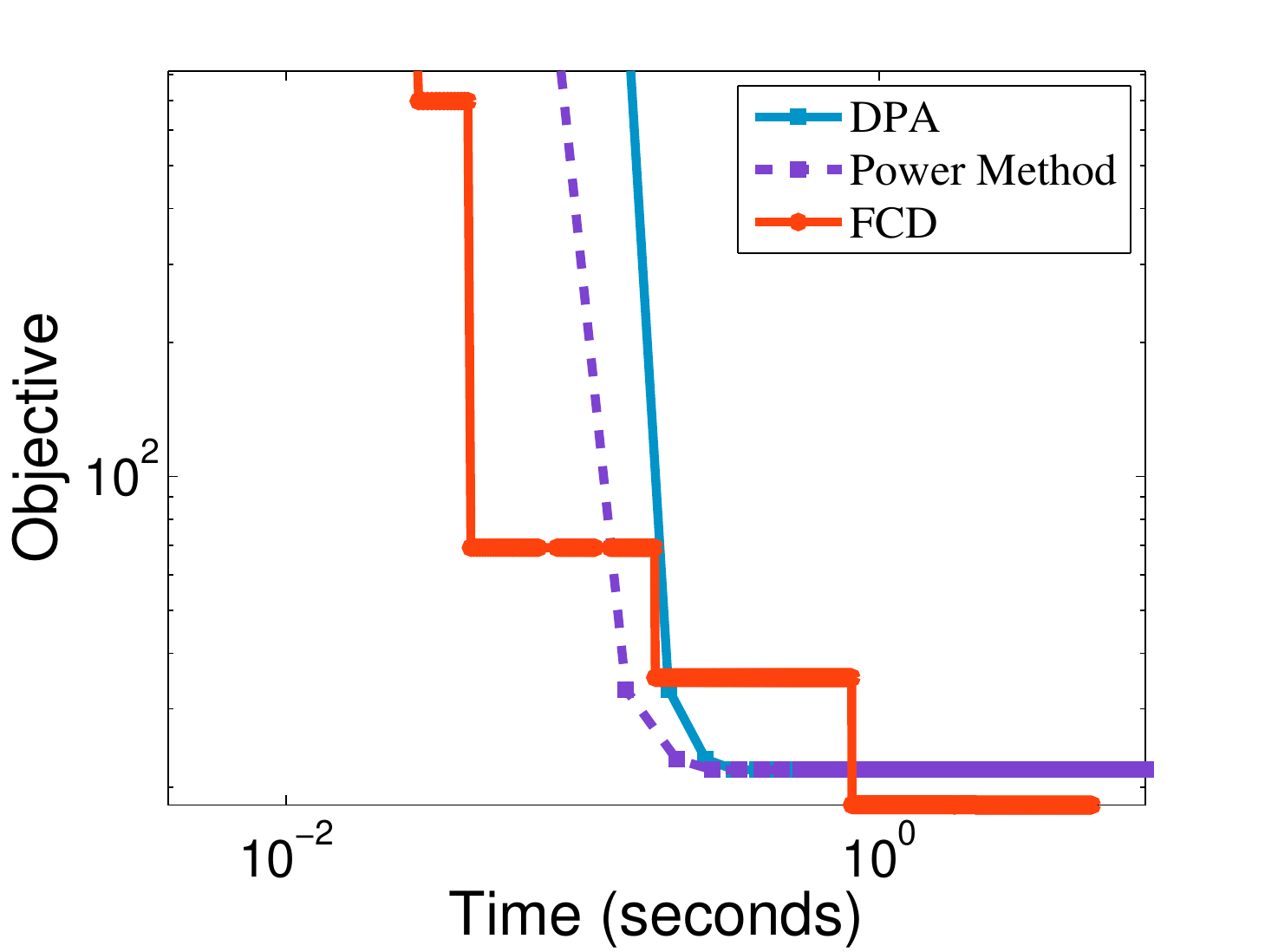}}
\caption{The convergence curve for solving the $\ell_p$ Norm Eigenvalue Problem with $p=4$.}\label{figure:ICA}
\end{figure}

\subsection{Experiment Results}

Table \ref{tab:acc:sparse:recovery} and Figure \ref{figure:sparse} show the accuracy and computational efficiency for the sparse recovery problem with setting $k=100$ and $\gamma = 0.1 / m$. We make the following observations. \textbf{(i)} The proposed method \textit{\textbf{PCD}} converges faster than the other methods. \textbf{(ii)} \textit{\textbf{PCD}} consistently gives the best performance.


Table \ref{tab:acc:ICA} and Figure \ref{figure:ICA} show the accuracy and computational efficiency for the $\ell_p$ Norm Eigenvalue Problem with $p=4$. We make the following observations. \textbf{(i)} Both PGSA and the power method present the same accuracy since they are essentially equivalent. \textbf{(ii)} While the other methods get stuck into poor local minima, \textit{\textbf{FCD}} exploits possible higher-order information of the non-convex function to escape from poor local minima and consistently finds lower objectives. This is consistent with our theory that our methods find stronger stationary points.

%

\clearpage
\normalem

\bibliographystyle{icml2022}
\bibliography{my}

\clearpage
\onecolumn
{
\LARGE
Appendix

}

The appendix is organized as follows. \\
Appendix \ref{app:sebsect:technical} presents the proofs for Section \ref{sebsect:technical}.\\
Appendix \ref{app:sebsect:convex} presents the proofs for Section \ref{sebsect:convex}.\\
Appendix \ref{app:sebsect:concave} presents the proofs for Section \ref{sebsect:concave}.\\
Appendix \ref{app:sebsect:discussion} presents some additional discussions.\\

\appendix

\section{Proofs for Section \ref{sebsect:technical}} \label{app:sebsect:technical}

\subsection{Proof of Lemma \ref{lemma:suff:dec}}

\begin{proof}

Noticing that $\nabla f(\cdot)$ is coordinate-wise Lipschitz continuous, we have:
\beq
f(\x^{t+1}) \leq f(\x^t) +  \la  \nabla  f(\x^t) ,\x^{t+1} - \x^t\ra  + \frac{1}{2} \|\x^{t+1} - \x^t\|_{\c}^2. \label{eq:lip:f}
\eeq

\textbf{(a)} We first discuss the \textit{\textbf{FCD}} algorithm. Using the optimality condition of (\ref{eq:subprob:nonconvex}), we have:
\beq
\frac{ f(\x^t) + \la \nabla f(\x^t),~\x^{t+1} - \x^t \ra  + \frac{\c_i+\theta}{2} \|\x^{t+1} - \x^t\|_2^2 + h(\x^{t+1})}{g(\x^{t+1})}\leq \frac{f(\x^t) + h(\x^t)}{g(\x^t)}.\nn
\eeq

\noi Combining this inequality with (\ref{eq:lip:f}), we have:
\beq
F(\x^{t+1} )+  \frac{\theta \|\x^{t+1}-\x^{t}\|_2^2 }{ 2 g(\x^{t+1})} \leq F(\x^t) .\nn
\eeq

\textbf{(b)} We now discuss the \textit{\textbf{PCD}} algorithm. Using the optimality condition of (\ref{eq:subprob:nonconvex2}) and the fact that $f(\x^{t}) + h(\x^{t}) = F(\x^{t}) \cdot g(\x^{t})$, we have:
\beq
&&f(\x^t) + \la  \nabla f(\x^t),\x^{t+1}-\x^{t} \ra + \frac{\c_i\eta^2}{2}  +   h(\x^{t+1}) + \frac{\theta}{2} \eta^2    - F(\x^t) g(\x^{t+1})  \nn\\
&\leq&
f(\x^t)   +   h(\x^t )   - F(\x^t) g(\x^t) =   0. \nn
\eeq

Rearranging terms, we have:
\beq
0  &\leq& - f(\x^t) - \la  \nabla f(\x^t),\x^{t+1}-\x^{t} \ra - \frac{\c_i+\theta}{2}\|\x^{t+1}-\x^{t}\|_2^2  -   h(\x^{t+1})    + F(\x^t) g(\x^{t+1})  \nn\\
   &\overset{(a)}{ \leq } &  - f(\x^{t+1}) - h(\x^{t+1}) - \frac{\theta}{2} \|\x^{t+1}-\x^{t}\|_2^2   + F(\x^t) g(\x^{t+1}) \nn\\
   &\overset{(b)}{ = } &   - F(\x^{t+1}) g(\x^{t+1})    + F(\x^t) g(\x^{t+1}) -  \frac{\theta}{2} \|\x^{t+1}-\x^{t}\|_2^2,\nn
\eeq
\noi where step $(a)$ uses (\ref{eq:lip:f}); step $(b)$ uses the definition $f(\x^{t+1}) + h(\x^{t+1}) = F(\x^{t+1}) \cdot g(\x^{t+1})$. Dividing both sides by $g(\x^{t+1})$ with $g(\x^{t+1})>0$, we have: $F(\x^{t+1})  - F(\x^{t}) \leq   - \frac{\theta}{2 g(\x^{t+1})} \|\x^{t+1} - \x^t\|^2$.

\end{proof}

\subsection{Proof of Lemma \ref{lemma:omega:sandwich}}

\begin{proof}
We now give a upper bound for $\alpha^t$. We have:
\beq
\alpha^t  &\triangleq  & \frac{ h(\x^{t+1}) +  f(\x^t) +  \la  \nabla f(\x^t),\x^{t+1} - \x^t\ra + \frac{\c_{i^t}+\theta}{2} \|\x^{t+1} - \x^t\|_2^2 }{g(\x^{t+1})} \nn\\
&\overset{(a)}{ \leq }&\frac{ h(\x^{t+1}) +   f(\x^{t+1}) + \frac{\c_{i^t}+\theta}{2} \|\x^{t+1} - \x^t\|_2^2  }{g(\x^{t+1})} \nn\\
& \overset{(b)}{ \leq } & F(\x^{t+1} ) +  ( F^{t} - F^{t+1} ) \cdot \frac{(\c_{i^t}+\theta)}{\theta},\nn
\eeq
\noi where step $(a)$ uses the convexity of $f(\cdot)$; step $(b)$ uses the sufficient decrease condition that: $\frac{\theta}{2 g(\x^{t+1})} \|\x^{t+1} - \x^t\|_2^2  \leq  F(\x^{t}) -F(\x^{t+1})$ as shown in Lemma \ref{lemma:suff:dec}.

We now give a lower bound for $\alpha^t$. We have:
\beq
\alpha^t  &\triangleq  & \frac{ h(\x^{t+1}) +  f(\x^t) +  \la  \nabla f(\x^t),\x^{t+1} - \x^t\ra + \frac{\c_{i^t}+\theta}{2} \|\x^{t+1} - \x^t\|_2^2 }{g(\x^{t+1})}\nn\\
&\overset{(a)}{ \geq }& \frac{h(\x^{t+1}) + f(\x^{t+1})  +  \frac{\theta}{2} \|\x^{t+1} - \x^t\|_2^2}{g(\x^{t+1})} \nn\\
&=& F(\x^{t+1}) + \frac{\theta \|\x^{t+1} - \x^t\|_2^2 }{2 g(\x^{t+1})} \geq F(\x^{t+1}),\nn
\eeq
\noi where step $(a)$ uses the fact that $\nabla f(\x)$ is coordinate-wise Lipschitz continuous.
\end{proof}

\subsection{Proof of Lemma \ref{label:property}}
\begin{proof}

First, since $\tilde{g}(\x)=-g(\x)$ is locally $\rho$-bounded non-convex at the point $\ddot{\x}$, applying Assumption \ref{ass:1} with $\x=\ddot{\x}$ and $\y=\ddot{\x}+\etas$ we have:
\beq\label{eq:weakly:convex:g}
&&-g(\ddot{\x})   \leq -g(\ddot{\x}+\etas)  - \la \ddot{\x} - (\ddot{\x}+\etas),~ \partial g(\ddot{\x}) \ra + \frac{\rho}{2}\|(\ddot{\x}+\etas) - \ddot{\x}\|_2^2 \nn\\
&\Rightarrow & \la \etas,~ \partial g(\ddot{\x}) \ra  \geq g(\ddot{\x}+\etas)-g(\ddot{\x}) - \frac{\rho}{2}\|\etas\|_2^2.
\eeq

Second, we have the following inequalities:
\begin{eqnarray} \label{eq:prop:gx}
 \forall  \etas,\sum_{i=1}^n g(\ddot{\x} + \etas_i \ei) & \overset{(a)}{ \geq }& ~\sum_{i=1}^n [ g(\ddot{\x}) + \la \partial g(\ddot{\x}),~ \etas_i \ei\ra]  \nn \\
 &\overset{(b)}{ = }& ~ n g(\ddot{\x}) + \la \partial g(\ddot{\x}),\etas\ra  \nn\\
 &\overset{(c)}{ \geq }& ~ n g(\ddot{\x}) +  g(\ddot{\x}+\etas)-g(\ddot{\x}) - \frac{\rho}{2}\|\etas\|_2^2,
\end{eqnarray}
\noi where step $(a)$ uses the convexity of $g(\cdot)$ that: $g(\x) - g({\x} + \etas_i \ei) +  \la ({\x} + \etas_i \ei)-\x,\partial g(\x) \ra \leq 0$; step $(b)$ uses $\la \partial g({\x}),~ \etas_i \ei\ra = \la\partial g({\x}), \etas\ra $; step $(c)$ uses (\ref{eq:weakly:convex:g}).

Third, we obtain the following equalities:
\beq\label{eq:prop:hx}
 \forall  \etas,\sum_{i=1}^n h(\ddot{\x} + \etas_i \ei) & = & \sum_{i=1}^n \left( h_i(\ddot{\x}_i + \etas_i) + \sum_{j\neq i} h_j(\ddot{\x}_j) \right) \nn\\
&=& \sum_{i=1}^n \left( h_i(\ddot{\x}_i + \etas_i) \right) + \sum_{i=1}^n \sum_{j\neq i} h_j(\ddot{\x}_j) \nn\\
&=& h(\ddot{\x}+\etas) + (n-1) h(\ddot{\x}).
\eeq


\textbf{(a)}  Since $\ddot{\x}$ is a \textit{FCW}-point, for all $\etas_i\in \mathbb{R}$, we have:
\beq
\mathcal{K}_i(\ddot{\x},0) &\leq& \mathcal{K}_i(\ddot{\x},\etas_i)   \nn\\
 \Leftrightarrow ~~  F(\ddot{\x}) &\leq& \frac{f(\ddot{\x}) +  \la  \nabla f(\ddot{\x}),~ \etas_i \ei \ra + h(\ddot{\x} +  \etas_i \ei) + \frac{\c_{i} + \theta}{2} \etas_i^2  }{g(\ddot{\x} + \etas_i \ei)}\nn\\
 \Leftrightarrow ~~ F(\ddot{\x})  g(\ddot{\x} + \etas_i \ei) &\leq& f(\ddot{\x}) +  \la  \nabla f(\ddot{\x}),~ \etas_i \ei \ra + h(\ddot{\x} + \etas_i \ei) + \frac{\c_{i} + \theta}{2} \etas_i^2.  \label{eq:same:form}
\eeq

Summing the inequality in (\ref{eq:same:form}) over $i = 1,...,n$, we have:
\beq \label{eq:sum:F}
\sum_{i=1}^n F(\ddot{\x}) \cdot g(\ddot{\x} + \etas_i \ei) \leq n f(\ddot{\x}) + \sum_{i=1}^n \la  \nabla f(\ddot{\x}),~ \etas_i \ei \ra + \sum_{i=1}^n h(\ddot{\x} + \etas_i \ei) +\sum_{i=1}^n \frac{\c_{i} + \theta}{2} \etas_i^2.
\eeq

\noi Combing (\ref{eq:prop:gx}), (\ref{eq:prop:hx}), and (\ref{eq:sum:F}), we have
\beq
&& F(\ddot{\x}) [n g(\ddot{\x}) +  g(\ddot{\x}+\etas)-g(\ddot{\x}) - \frac{\rho}{2}\|\etas\|_2^2]\nn\\
 &\leq&  n f(\ddot{\x}) + \sum_{i=1}^n \la  \nabla f(\ddot{\x}), \etas_i \ei \ra + \sum_{i=1}^n \frac{\c_{i} + \theta}{2} \etas_i^2 +  h(\ddot{\x}+\etas) + (n-1) h(\ddot{\x}) \nn\\
 &\overset{(a)}{ = }&n f(\ddot{\x}) +  \la  \nabla f(\ddot{\x}),  \etas    \ra + \frac{1}{2}\| \etas \|^2_{\c+\theta} +  h(\ddot{\x}+\etas) + (n-1) h(\ddot{\x})  \nn\\
 &\overset{(b)}{ \leq } &n f(\ddot{\x}) +  f(\ddot{\x}+\etas) - f(\ddot{\x}) + \frac{1}{2}\| \etas \|^2_{\c+\theta}  +  h(\ddot{\x}+\etas) + (n-1) h(\ddot{\x}) \nn\\
 &\overset{ }{ = } & (n-1) ( f(\ddot{\x}) + h(\ddot{\x})) +  f(\ddot{\x}+\etas)+  h(\ddot{\x}+\etas) + \frac{1}{2}\| \etas \|^2_{\c+\theta},     \label{fcd:finall}
\eeq
\noi where step $(a)$ uses $\sum_{i=1}^n \la  \nabla f(\ddot{\x}),~ \etas_i \ei \ra=\la  \nabla f(\ddot{\x}), \etas    \ra$; step $(b)$ uses the convexity of $f(\cdot)$ that:
\beq
f(\x) - f({\x} + \etas_i \ei) +  \la ({\x} + \etas_i \ei)-\x,\nabla  f(\x) \ra \leq 0. \nn
\eeq

Finally, from (\ref{fcd:finall}) we have the following results:
\beq
&& F(\ddot{\x}) \cdot [ n g(\ddot{\x}) +    g(\ddot{\x}+\etas) -    g(\ddot{\x}) ] \leq  (n-1)  F(\ddot{\x}) g(\ddot{\x})   +  f(\ddot{\x}+\etas)+  h(\ddot{\x}+\etas)    + \mathcal{C}(\ddot{\x},\etas) \nn\\
&\Leftrightarrow&    F(\ddot{\x})    g(\ddot{\x})  \leq   f(\ddot{\x}+\etas)    +   h(\ddot{\x}+\etas) + F(\ddot{\x})(g(\ddot{\x}) - g(\ddot{\x} + \etas))  + \mathcal{C}(\ddot{\x},\etas)  \nn\\
&\Leftrightarrow&    F(\ddot{\x})   g(\ddot{\x}+\etas)  \leq   f(\ddot{\x}+\etas) +h(\ddot{\x}+\etas)    + \mathcal{C}(\ddot{\x},\etas)  \nn\\
&\Leftrightarrow&    F(\ddot{\x})   \leq   F(\ddot{\x}+\etas)  + \frac{ \mathcal{C}(\ddot{\x},\etas)}{g(\ddot{\x}+\etas)} \label{eq:same:optimality}
\eeq

\noi Therefore, we finish the first part of this lemma.

\noi \textbf{(b)}  Since $\dot{\x}$ is a \textit{PCW}-point, for all $\etas_i\in \mathbb{R}$, we have:
\beq
\mathcal{M}_i(\dot{\x},0) &\leq& \mathcal{M}_i(\dot{\x},\etas_i)   \nn\\
 \Leftrightarrow ~~  f(\dot{\x}) + h(\dot{\x})   - F(\dot{\x}) \cdot g(\dot{\x}  ) &\leq& (\Q_{i}(\dot{\x},\eta) +  h(\dot{\x} + \etas_i \ei) + \frac{\theta}{2} \etas_i^2)   - F(\dot{\x}) \cdot g(\dot{\x} + \etas_i \ei) \nn\\
 \Leftrightarrow ~~ F(\dot{\x}) \cdot g(\dot{\x} + \etas_i \ei) &\leq& f(\dot{\x}) +  \la  \nabla f(\dot{\x}),~ \etas_i \ei \ra +  h(\dot{\x} + \etas_i \ei)+\frac{\c_{i} + \theta}{2} \etas_i^2 \nn
\eeq
The inequality above has the same form as in (\ref{eq:same:form}). Therefore, we have a similar conclusion to (\ref{eq:same:optimality}) that: $F(\dot{\x})\leq   F(\ddot{\x}+\etas) + \frac{\mathcal{C}(\dot{\x},\etas)}{g(\dot{\x}+\etas)}$.

\end{proof}
\section{Proofs for Section \ref{sebsect:convex}} \label{app:sebsect:convex}

\subsection{Proof of Theorem \ref{the:optimality}}

\begin{proof}

(\textbf{a}) $\{$Optimal point $\bar{\x}\}$ $\in$ $\{$\textit{FCW}-point $\ddot{\x}\}$. By the optimality of $\bar{\x}$, we have:
\beq
\frac{f(\bar{\x}) +   h(\bar{\x})}{   g(\bar{\x})}    \leq \frac{f(\x) +   h(\x)}{g(\x)},~\forall \x \nn
\eeq
\noi Letting $\x = \bar{\x} + \etas_i \ei$, we have:
\beq \label{eq:opt:to:fcw}
\frac{f(\bar{\x}) +   h(\bar{\x})}{   g(\bar{\x})} &\leq & \frac{f(\bar{\x} + \etas_i \ei ) +   h(\bar{\x} + \etas_i \ei)}{g(\bar{\x} + \etas_i \ei)},~\forall \etas_i \nn \\
&\overset{(a)}{ \leq } &  \frac{ f(\bar{\x}) +  \la  \nabla_i f(\bar{\x}),~ \etas_i \ei ) + \frac{\c_i}{2} \etas_i^2 +   h(\bar{\x} + \etas_i \ei)}{g(\bar{\x} + \etas_i \ei)},~\forall \etas_i \nn \\
&\overset{(b)}{ \leq } &  \frac{ f(\bar{\x}) +  \la  \nabla_i f(\bar{\x}),~ \etas_i \ei ) + \frac{\c_i}{2} \etas_i^2 + \frac{\theta}{2} \etas_i^2 +   h(\bar{\x} + \etas_i \ei)}{g(\bar{\x} + \etas_i \ei)},~\forall \etas_i \nn \\
&\overset{(c)}{ = } & \K_{i}(\bar{\x},\etas_i ) ,~\forall \etas_i,
\eeq
\noi where step $(a)$ uses coordinate-wise Lipschitz continuity of $\nabla f(\cdot)$ that: $f(\bar{\x} + \etas_i \ei  ) \leq f(\bar{\x}) +  \la  \nabla_i f(\bar{\x}),~ \etas_i \ei ) + \frac{\c_i}{2} \etas_i^2,~\forall \etas_i$; step $(b)$ uses the fact that $\theta>0$; step $(c)$ uses the definition of $\K_{i}(\bar{\x},\etas_i)$. Using the fact that $\mathcal{K}_i(\bar{\x},0) = \frac{f(\bar{\x}) +   h(\bar{\x})}{   g(\bar{\x})}$. The inequality in (\ref{eq:opt:to:fcw}) essentially implies that:
\beq
\mathcal{K}_i(\bar{\x},0) = \min_{\etas_i}\K_{i}(\bar{\x},\etas_i ).   \nn
\eeq
\noi Therefore, any optimal point $\bar{\x}$ must be a \textit{FCW}-point.

(\textbf{b}) $\{$\textit{FCW}-point $\dot{\x}\}$ $\Leftrightarrow$ $\{$\textit{PCW}-point $\ddot{\x}\}$. Using the optimality of \textit{FCW}-point and \textit{PCW}-point, we respectively have the following inequalities:
\beq
F(\ddot{\x}) \cdot g(\ddot{\x} + \etas_i \ei) &\leq& f(\ddot{\x}) +  \la  \nabla f(\ddot{\x}),~ \etas_i \ei \ra + h(\ddot{\x}+ \etas_i \ei ) + \frac{\c_{i} + \theta}{2} \etas_i^2 ,\forall \etas_i; \nn\\
F(\dot{\x}) \cdot g(\dot{\x} + \etas_i \ei) &\leq& f(\dot{\x}) +  \la  \nabla f(\dot{\x}),~ \etas_i \ei \ra + h(\dot{\x}+ \etas_i \ei ) + \frac{\c_{i} + \theta}{2} \etas_i^2,\forall \etas_i. \nn
\eeq
\noi These two inequalities have the same form, leading to the same optimality condition as shown in Lemma \ref{label:property}. We conclude that the optimality of \textit{FCW}-point is completely equivalent to that of \textit{PCW}-point.

(\textbf{c}) $\{$\textit{FCW}-point $\ddot{\x}$\} $\in$ $\{$\textit{D}-point $\grave{\x}\}$. By assumption, we have $g(\x)\geq \underline{g} > 0$ for all $\x$ for some universal constant $\underline{g}$. For any $\y \in\text{dom}(F)$, we let $\etas = t(\y -\ddot{\x} )$ and have the following results:
\beq
&& \lim_{t\downarrow 0} \frac{1}{t} \cdot \left( F(\ddot{\x} + t(\y -\ddot{\x} ) ) - F(\ddot{\x}) \right) \nn\\
&=& \lim_{t\downarrow 0} \frac{1}{t} \cdot \left( F(\ddot{\x}+\etas)  - F(\ddot{\x}) \right) \nn\\
&\overset{(a)}{ \geq }&\lim_{t\downarrow 0} -\frac{1}{t} \cdot \frac{\mathcal{C}(\ddot{\x},\etas)}{g(\ddot{\x}+\etas)}  \nn\\
&\overset{(b)}{ \geq }&\lim_{t\downarrow 0} -\frac{1}{t\underline{g}} \cdot \mathcal{C}(\ddot{\x},\etas) \nn\\
&\overset{(c)}{ \geq }&\lim_{t\downarrow 0} -\frac{1}{t\underline{g}} \cdot [ \frac{1}{2 }\| \etas \|^2_{(\c+\theta)} + \frac{ \rho}{2}\| \etas \|^2 \cdot  F(\ddot{\x})  ]\nn\\
&\overset{(d)}{ \geq }&\lim_{t\downarrow 0} -\frac{1}{t \underline{g}} \cdot [ \frac{t^2}{2 } \| \y -\ddot{\x}  \|^2_{(\c+\theta)} + \frac{ \rho t^2}{2}\| \y -\ddot{\x} \|^2 \cdot  F(\ddot{\x})  ]\nn\\
&=&0,\nn
\eeq
\noi where step $(a)$ uses the property of \textit{FCW}-point as in Lemma \ref{label:property}; step $(b)$ uses the assumption that $g(\ddot{\x}+\etas) \geq  \underline{g}$; step $(c)$ uses the definition of $\mathcal{C}(\x,\etas)  \triangleq  \frac{1}{2}\| \etas \|^2_{\c+\theta} +  \frac{\rho}{2}\|\etas\|_2^2  F(\x) $; step $(d)$ uses $\etas = t(\y -\ddot{\x} )$. Therefore, any \textit{FCW}-point must be a \textit{D}-point.

(\textbf{c}) $\{$\textit{D}-point $\grave{\x}$\} $\in$ $\{$\textit{C}-point $\check{\x}\}$. We define $\z \triangleq \grave{\x} + t(\y-\grave{\x})$ and derive the following inequalities:
\beq \label{eq:hhhh}
0 &\overset{(a)}{ \leq }& \lim_{t\downarrow 0}~\frac{1}{t} \cdot (F( \z) - F(\grave{\x}) ) \nn\\
&\overset{(b)}{ = }& \lim_{t\downarrow 0}~\frac{1}{tg(\z)}\cdot [ f(\z)+h(\z)  - F(\grave{\x}) \cdot g(\z)   ]  \nn\\
& \overset{(c)}{ \leq } & \lim_{t\downarrow 0}~ \frac{1}{tg(\z)}\cdot[    f(\grave{\x}) + h(\grave{\x}) + \la \z-\grave{\x},~ \nabla f(\z) + \partial h(\z)\ra      + (-g(\grave{\x}) + \la \grave{\x}-\z, \partial g(\grave{\x}) \ra )\cdot F(\grave{\x})]  \nn\\
& \overset{(d)}{ = } & \lim_{t\downarrow 0}~ \frac{1}{tg(\z)}\cdot[   \la \z-\grave{\x},~ \nabla f(\z) + \partial h(\z)\ra      + \la \grave{\x}-\z, \partial g(\grave{\x}) \ra F(\grave{\x})]  \nn\\
& \overset{(e)}{ = } & \lim_{t\downarrow 0}~ \frac{1}{tg(\z)}\cdot   \la t(\y-\grave{\x}),~ \nabla f(\z) + \partial h(\z)   - F(\grave{\x}) \partial g(\grave{\x}) \ra   \nn\\
& \overset{(f)}{ = } & \lim_{t\downarrow 0}~ \frac{1}{g(\grave{\x})}\cdot   \la \y-\grave{\x},~ \nabla f(\grave{\x}) + \partial h(\grave{\x})   - F(\grave{\x}) \partial g(\grave{\x}) \ra,
\eeq
\noi where step $(a)$ uses the definition of \textit{D}-point that: $0 \leq \lim_{t\downarrow 0}~\frac{1}{t}[F( \grave{\x} + t(\y-\grave{\x})) - F(\grave{\x}) ]$; step $(b)$ uses the definition of $F(\x) = \frac{f(\x)+h(\x)}{g(\x)}$; step $(c)$ uses the convexity of $f(\cdot), h(\cdot)$ and $g(\cdot)$ that:
\beq
f(\z) \leq f(\grave{\x}) + \la \z-\grave{\x},~\nabla f(\z) \ra ,\nn\\
h(\z) \leq h(\grave{\x}) + \la \z-\grave{\x},~ \partial h(\z) \ra ;\nn \\
-g(\z) \leq  - g(\grave{\x}) + \la \grave{\x}-\z,~ \partial g(\grave{\x}) \ra ;\nn
\eeq
\noi step $(d)$ uses the definition of $F(\x) = \frac{f(\x)+h(\x)}{g(\x)}$; step $(e)$ uses $\z - \grave{\x} = t(\y-\grave{\x})$; step $(f)$ uses $\z=\grave{\x}$ with $t\downarrow0$;

Noticing that $g(\grave{\x})>0$ amd the inequality in (\ref{eq:hhhh}) holds for all $\y \in \text{dom}(F)$, we have:
\beq
0 \in \nabla f(\grave{\x}) + \partial h(\grave{\x}) - F(\grave{\x})\cdot \partial g(\grave{\x}).\nn
\eeq
\noi Therefore, any \textit{D}-point must be a \textit{C}-point.

\end{proof}

\subsection{Proof of Proposition \ref{prop:convex:conv}}

\begin{proof}

First, we note that the sequence $\{F(\x^t)\}_{t\geq0}$ is monotonically non-increasing. Taking the expectation for both sides of the sufficient decrease condition as shown in Lemma \ref{lemma:suff:dec}, we have:
\beq
\E_{i^t}[F^{t+1}] - F(\x^t) \leq -\frac{\theta}{n g(\x^{t+1})}\E_{i^t}[\|\x^{t+1} - \x^t\|_2^2]. \nn
\eeq
\noi Summing the inequality above over $t=0,1,...,T$, we have:
\beq
\E_{\xi^{T}} [\sum_{t=0}^T \frac{ \theta \|\x^{t+1}-\x^{t}\|_2^2 }{2 g(\x^{t+1})}] \leq \E_{\xi^T} [n(F(\x^0) - F(\x^{T+1}))] \leq n(F(\x^0) - F(\bar{\x})) . \label{eq:global:conv:convex}
\eeq

\noi Combining with the fact that $\g(\x^t)\leq \bar{g}$ and $F(\bar{\x})\geq0$, we conclude that
\beq
\E_{\xi^{T}}[ \sum_{t=0}^T \| \x^{i+1}-\x^{i} \|_2^2 ] \leq \tfrac{2 n \bar{g} F(\x^0) }{\theta(T+1)}.\nn
\eeq

Therefore, there exists an index $\bar{t}$ with $0\leq \bar{t}\leq T$ such that:
\beq\label{eq:above:1}
\E_{\xi^{T}} [ \| \x^{\bar{t}+1}- \x^{\bar{t}}\|_2^2 ] \leq   \frac{2 n \bar{g} F(\x^0)  }{\theta (T+1)}
\eeq

\noi We have $\lim_{t\rightarrow \infty}\E_{\xi^{t}}[\|\x^{t+1} - \x^t\|_2^2] = 0$. Therefore, every clustering point of the sequence of \textit{\textbf{FCD}} is almost surely a \textit{FCW}-point of Problem (\ref{eq:main}).

Furthermore, for any $\bar{t}$, we have:
\beq \label{eq:above:2}
\E_{\xi^{\bar{t}}} [ \| \x^{\bar{t}+1}- \x^{\bar{t}}\|_2^2 ]  = \frac{1}{n}\sum_{i=1}^n \text{dist}(0,\arg \min_{\eta} \mathcal{K}_{i}(\x^{\bar{t}},\eta)   )^2
\eeq

Combining (\ref{eq:above:2}) and (\ref{eq:above:1}), we have the following result:
\beq
\frac{1}{n}\sum_{i=1}^n \text{dist}(0,\arg \min_{\eta} \mathcal{K}_{i}(\x^{\bar{t}},\eta)   )^2 \leq  \frac{2 n \bar{g} F(\x^0)  }{\theta (T+1)}\nn
\eeq

\noi We conclude that \textit{\textbf{FCD}} finds an $\epsilon$-approximate \bfit{FCD}-point in at most $T+1$ iterations in the
sense of expectation, where
\beq
T \leq \lceil \frac{2 n \bar{g} F(\x^0)  }{\theta \epsilon } \rceil = \mathcal{O}(\epsilon^{-1}).\nn
\eeq


Using similar strategy, we can prove that \textit{\textbf{PCD}} converges to a \textit{\textbf{PCW}}-point whenever \textit{\textbf{PCD}} converges.

\end{proof}

\subsection{Proof of Theorem \ref{the:rate:FCDC}}

\begin{proof}
We prove the convergence rate of \textit{\textbf{FCD}} for convex-convex FMPs.

We define
\beq \label{eq:definition:varpi}
\bar{\rho} \triangleq \frac{\rho}{\min(\bar{\c})}, \varpi \triangleq \left(\frac{\rho}{\min(\bar{\c})}\right)\cdot \left(\frac{\max(\bar{\c})}{\theta} F(\x^0)\right).
\eeq

\noi First, for any $\x,\d\in \mathbb{R}^n$, we have the following equalities:
\beq
\frac{1}{n}\sum_{i=1}^n\|\x + \d_i \ei\|_2^2 &=& \frac{1}{n}\|\d\|_2^2 + \frac{2}{n} \la \x, \d \ra +  \|\x\|_2^2 \nn\\
&=& \frac{1}{n}\|\d + \x\|_2^2  + (1-\frac{1}{n})\|\x\|_2^2 \nn
\eeq
\noi Applying the equality above with $\x=\x^{t} - \bar{\x}$ and $\d=\x^{t+1}-\x^t$, we have:
\beq \label{eq:imp:exp}
\E_{i^{t}}[ \|\x^{t+1} - \bar{\x}\|_{\bar{\c}}^2] =  \frac{1}{n}\| \x^{t+1} - \bar{\x}\|_{\bar{\c}}^2 + (1-\frac{1}{n}) (\|\x^{t} - \bar{\x}\|_{\bar{\c}}^2).
\eeq

\noi Second, the optimality condition for the non-convex subproblem as in (\ref{eq:subprob:nonconvex}) can be written as:
\beq
  &&0 \in    [ \nabla_{i^t} f(\x^t) + \bar{\c}_{i^t}  \bar{\eta}^t  + \partial_{i^t} h(\x^{t+1}) ]  g(\x^{t+1})  - \J_{i^t}(\x^t,\bar{\eta}^t,\theta) \cdot \partial_{i^t} g(\x^{t+1})  \nn \\
  &\Leftrightarrow & 0 \in  \nabla_{i^t} f(\x^t) + \bar{\c}_{i^t} \bar{\eta}^t + \partial_{i^t} h(\x^{t+1})  - \alpha^t \partial_{i^t} g(\x^{t+1}) . \label{eq:fractional:cd:convex:optimality:condition}
\eeq

\noi For any $\x \in\mathbb{R}^n$, we derive the following results:
\beq \label{eq:fractional:cd:convex:0}
&&\E_{i^t}[\frac{1}{2}\|\x^{t+1}-\x\|_{\bar{\c}}^2] -  \E\frac{1}{2}\|\x^{t}-\x\|_{\bar{\c}}^2]\nn\\
&\overset{(a)}{ = }& \E_{i^{t}}[\la \x-\x^{t+1},\bar{\c} \odot (\x^{t}- \x^{t+1})  \ra ] - \E_{i^{t}}[\frac{1}{2} \la \x^{t}- \x^{t+1}, \bar{\c} \odot (\x^{t}- \x^{t+1}) \ra ] \nn\\
&\overset{(b)}{ = }&  \E_{i^{t}}[\la  \x-\x^{t+1},          (\nabla_{i^t} f(\x^t) + \partial_{i^t} h(\x^{t+1}) - \alpha^t \partial_{i^t} g(\x^{t+1}) ) \cdot \eit  \ra ] \nn\\
 && - \E_{i^{t}}[\frac{1}{2} \la \x^{t}- \x^{t+1}, (\bar{\c}_{i^t}  (\x^{t}_{i^t}- \x^{t+1}_{i^t})) \cdot \eit \ra ]  \nn\\
&\overset{(c)}{ = }&    \frac{1}{n} \la \x - \x^{t+1},\nabla f(\x^t) + \partial h(\x^{t+1}) - \alpha^t \partial g(\x^{t+1})     \ra  + \frac{1}{2n} \| \x^{t+1}-\x^{t} \|_{\bar{\c}}^2 \nn \\
&\overset{}{ = }&   \frac{1}{n} [\la \x-\x^{t+1},\nabla f(\x^t) \ra +  \la   \x-\x^{t+1} ,\partial h(\x^{t+1}) - \alpha^t\partial g(\x^{t+1})  \ra ] + \frac{1}{2n} \| \x^{t+1}-\x^{t} \|_{\bar{\c}}^2 \nn\\
&\overset{(d)}{ \leq }&    \frac{1}{n} [ \la \x-\x^{t+1},\nabla f(\x^t) + \partial h(\x^{t+1}) \ra ]   \nn\\
  && + \frac{\alpha^t}{n}[  g(\x^{t+1})- g(\x)  +   \frac{\bar{\rho}}{2}\|\x-\x^{t+1}\|_{\bar{\c}}^2  ]+ \frac{1}{2n} \| \x^{t+1}-\x^{t} \|_{\bar{\c}}^2,
\eeq
\noi where step $(a)$ uses the Pythagoras relation that: $\forall \x,\y, \z,\frac{1}{2}\|\y-\z\|_2^2 - \frac{1}{2}\|\x-\z\|_2^2=  \la \z-\y,\x - \y\ra - \frac{1}{2}\|\x-\y\|_2^2$; step $(b)$ uses the optimality condition in (\ref{eq:fractional:cd:convex:optimality:condition}); step $(c)$ uses the fact that $\E_{i^{t}}[ \la \x_{i^t} \eit, \y \ra] = \frac{1}{n} \sum_{j=1}^n  \x_j \y_j = \frac{1}{n}\la \x,\y\ra$; step $(d)$ uses the convexity of $f(\cdot)$ that:
\beq
\la   \x - \x^t ,\nabla f(\x^t) \ra &\leq& f(\x)    - f(\x^{t}); \nn
\eeq
\noi step $(e)$ uses the $\rho$-bounded non-convexity of $-g(\cdot)$ that:
\beq\label{eq:fractional:cd:convex:weak:p}
-\la  \x- \x^{t+1} ,\partial g(\x^{t+1}) \ra   & \leq & - g(\x) + g(\x^{t+1}) +   \frac{\rho}{2}\|\x-\x^{t+1}\|_2^2 \nn \\
 & \leq & - g(\x) + g(\x^{t+1}) +   \frac{\rho   }{2 \min(\bar{\c})}\|\x-\x^{t+1}\|_{\bar{\c}}^2 \nn  \\
 & = & - g(\x) + g(\x^{t+1}) +   \frac{\bar{\rho}}{2}\|\x-\x^{t+1}\|_{\bar{\c}}^2 .\nn
\eeq

We further derive the following results:
\beq\label{eq:fractional:cd:convex:1}
&&\la \x-\x^{t+1},\nabla f(\x^t) + \partial h(\x^{t+1}) \ra + \alpha^t g(\x^{t+1}) \nn\\
& = & \la \x-\x^{t+1},\partial h(\x^{t+1}) \ra  + \la \x-\x^{t},\nabla f(\x^t) \ra  + \la \x^{t}-\x^{t+1},\nabla f(\x^t) \ra + \alpha^t g(\x^{t+1})  \nn\\
& \overset{(a)}{ \leq } & h(\x) - h(\x^{t+1})  + f(\x) - f(\x^{t})   + \la \x^{t}-\x^{t+1},\nabla f(\x^t) \ra + \alpha^t g(\x^{t+1})  \nn\\
& \overset{(b)}{ = } & h(\x) - h(\x^{t+1})  + f(\x) - f(\x^{t})   + \la \x^{t}-\x^{t+1},\nabla f(\x^t) \ra + \J_{i^t}(\x^t,\bar{\eta}^t,\theta)  \nn\\
& \overset{(c)}{ = } & h(\x) - h(\x^{t+1})  + f(\x) - f(\x^{t})   +   f(\x^t) +   h(\x^{t+1}) + \frac{\c_{i^t} + \theta}{2} \|\x^{t+1}-\x^{t}\|_2^2  \nn\\
& \overset{}{ = } & h(\x)   + f(\x)    + \frac{1}{2} \|\x^{t+1}-\x^{t}\|_{\bar{\c}}^2,
\eeq
\noi where step $(a)$ uses the convexity of $f(\cdot)$ and $h(\cdot)$; step $(b)$ uses the fact that $\alpha^t = \J_{i^t}(\x^t,\bar{\eta}^t,\theta) /g(\x^{t+1})$; step $(c)$ uses the definition of $\J_{i^t}(\x^t,\bar{\eta}^t,\theta)$ that $\J_{i^t}(\x^t,\bar{\eta}^t,\theta) = f(\x^t) + h(\x^{t+1}) + \la \x^{t+1}-\x^t,\nabla f(\x^t) \ra + \frac{\theta + \c_{i^t}}{2}\|\x^{t+1}-\x^t\|_2^2$.

Combining (\ref{eq:fractional:cd:convex:0}) and (\ref{eq:fractional:cd:convex:1}), we have:
\beq \label{eq:fractional:cd:convex:2}
&& \E_{i^{t}}[\frac{1}{2}\|\x^{t+1}-\x\|_{\bar{\c}}^2] - \E [ \frac{1}{2}\|\x^{t}-\x\|_{\bar{\c}}^2      ]  \nn\\
&\leq & \frac{1}{n} [h(\x)   + f(\x)    - \alpha^t g(\x)] +  \frac{\alpha^t \bar{\rho}}{2 n }\|\x-\x^{t+1}\|_{\bar{\c}}^2 \nn\\
&\overset{(a)}{ = } & \frac{g(\x)}{n} [F(\x)    - \alpha^t] +  \frac{\alpha^t \bar{\rho}}{2 n }\|\x-\x^{t+1}\|_{\bar{\c}}^2 \nn\\
& \overset{(b)}{ \leq } & \frac{g(\x)}{n} [F(\x)    - F(\x^{t+1})] +  \frac{\varpi}{2 n }\|\x-\x^{t+1}\|_{\bar{\c}}^2  \nn\\
& \overset{(c)}{ = } & \frac{g(\x)}{n} [F(\x)    - F(\x^{t+1})] +  \frac{\varpi}{2} \E_{i^{t}}[\|\x^{t+1}-\x\|_{\bar{\c}}^2] - (1-\frac{1}{n})\frac{\varpi}{2}\|\x^{t}-\x\|_{\bar{\c}}^2,
\eeq
\noi where step $(a)$ uses the fact that $g(\x)F(\x) = h(\x) + f(\x)$; step $(b)$ uses the inequality $F(\x^{t+1})\leq\alpha^t$ and $\alpha^t \bar{\rho} \leq \sigma F(\x^0) \cdot\bar{\rho}   \triangleq \varpi$ as shown in Lemma \ref{lemma:omega:sandwich} and (\ref{eq:definition:varpi}); step $(c)$ uses (\ref{eq:imp:exp}).

\noi We apply (\ref{eq:fractional:cd:convex:2}) with $\x=\ddot{\x}$ and rearranging terms, we obtain:
\beq
(1- \varpi) \E_{i^{t}}[r^{t+1}]  + \frac{g(\x)}{n}\ddot{q}^{t+1} &\leq& (1-\varpi) r^t + \frac{\varpi}{n} r^t    \label{eq:fractional:cd:convex:3}
\eeq

We now discuss the case when $F(\cdot)$ satisfies the Luo-Tseng error bound assumption. We first bound the term $r^t$ in (\ref{eq:fractional:cd:convex:3}) using the following inequalities:
\beq
r^t\nn & \triangleq &  \max(\bar{\c}) \frac{1}{2}\|\x^t - \ddot{\x}\|_{2}^2 \nn\\
&\overset{(a)}{ \leq } & \max(\bar{\c}) \frac{1}{2}\frac{\delta^2}{n^2} (   \sum_{i=1}^n |\P_i(\x^t)|)^2 \nn\\
&\overset{(b)}{ \leq } & \max(\bar{\c}) \frac{1}{2}\frac{\delta^2}{n^2} n \cdot (   \sum_{i=1}^n |\P_i(\x^t)|^2) \nn\\
&\overset{(c)}{ \leq } & \max(\bar{\c}) \frac{1}{2}\frac{\delta^2}{n^2} n \cdot \left(n \E_{i^t}[\|\x^{t+1}-\x^t\|^2_2]\right)  \nn\\
&\overset{(d)}{ = } &  \max(\bar{\c})  \frac{\delta^2 g(\x^{t+1})}{\theta} \cdot \frac{\theta}{2g(\x^{t+1})} \E[\|\x^{t+1}-\x^t\|^2_2] \nn\\
&\overset{(e)}{ \leq } &  \max(\bar{\c})\delta^2 \frac{\bar{g}}{\theta}(F(\x^t)-F(\x^{t+1})) \nn\\
& = &  \max(\bar{\c})\delta^2 \frac{\bar{g}}{\theta}(\ddot{q}^t-\ddot{q}^{t+1}),\label{eq:conv:rate:tsengluo:1}
\eeq
\noi where step $(a)$ uses the Luo-Tseng error bound assumption as in (\ref{eq:luotseng}); step $(b)$ uses the fact that $\|\x\|_1^2\leq n \|\x\|_2^2,~\forall \x \in \mathbb{R}^n$; step $(c)$ uses the fact that $\E_{i^t}[\|\x^{t+1}-\x^t\|^2_2] = \E_{i^t}[\|(\x^t + \P_{i^t}(\x^t))-\x^t\|^2_2] =  \E_{i^t}[|\P_{i^t}(\x^t)|^2] = \frac{1}{n}\sum_{i=1}^n |\P_i(\x^t)|^2$; step $(d)$ uses the assumption that $g(\x^{t+1})\leq \bar{g}$ and the sufficient decrease condition in Lemma \ref{lemma:suff:dec}.

Since $\varpi\leq 1$, we have form (\ref{eq:fractional:cd:convex:3}):
\beq
\frac{g(\bar{\x})}{n}\ddot{q}^{t+1} &\leq& (1-\varpi) r^t + \frac{\varpi}{n} r^t    \nn\\
&\overset{(a)}{ \leq } &  (1 + \frac{1}{n}) r^t  \nn\\
&\overset{(b)}{ \leq } &  (1 + \frac{1}{n}) \max(\bar{\c})\delta^2 \frac{\bar{g}}{\theta}(\ddot{q}^t-\ddot{q}^{t+1})  \nn\\
&\overset{(c)}{ = } & \frac{\kappa_1 \bar{g}}{n} (\ddot{q}^t-\ddot{q}^{t+1}) \label{eq:conv:rate:tsengluo:11}
\eeq
\noi where step $(a)$ uses $0<\varpi\leq 1$; step $(b)$ uses (\ref{eq:conv:rate:tsengluo:1}); step $(c)$ uses the definition of $\kappa_1$ that $\kappa_1  \triangleq (n + 1) \max(\bar{\c})\delta^2 \frac{1}{\theta}$.

Finally, using the definition of $\kappa_0$ that $\kappa_0 \triangleq \frac{g(\bar{\x})}{\bar{g}}$, we have the following results from (\ref{eq:conv:rate:tsengluo:11}):
\beq
&&\kappa_0 \ddot{q}^{t+1} \leq \kappa_1 (\ddot{q}^t -\ddot{q}^{t+1}) \nn\\
& \Rightarrow & \ddot{q}^{t+1} \leq \frac{\kappa_1}{\kappa_1 + \kappa_0} \ddot{q}^t \nn\\
& \Rightarrow & \ddot{q}^{t+1} \leq (\frac{\kappa_1}{\kappa_1 + \kappa_0})^{t+1} \ddot{q}^0 \nn
\eeq
\noi Thus, we finish the proof of this theorem.

\end{proof}

\subsection{Proof of Theorem \ref{the:rate:DCDC}}

\begin{proof}

We prove the convergence rate of \textit{\textbf{PCD}} for convex-convex FMPs.

We define $\bar{\rho} = \frac{\rho}{\min(\bar{\c})}$.

\noi The optimality condition for the non-convex subproblem as in (\ref{eq:subprob:nonconvex}) can be written as:
\beq \label{eq:pcd:convex:optimality:condition}
 0 \in   \nabla_{i^t} f(\x^t)  + \partial_{i^t} h(\x^{t+1}) + (\c_{i^t}+\theta) \bar{\eta}^t   - F(\x^t) \cdot  \partial_{i^t} g(\x^{t+1}).
\eeq

\noi \noi For any $\x \in\mathbb{R}^n$, we derive the following results:
\beq \label{eq:parametric:cd:convex:0}
&&\E_{i^{t}}[\frac{1}{2}\|\x^{t+1}-\x\|_{\bar{\c}}^2\ -  \E[\frac{1}{2}\|\x^{t}-\x\|_{\bar{\c}}^2]\nn\\
&\overset{(a)}{ = }& \E_{i^{t}}[\la \x-\x^{t+1},\bar{\c} \odot (\x^t - \x^{t+1})  \ra ] - \E_{i^{t}}[  \frac{1}{2}\|(\x^{t+1}- \x^t)\|_{\bar{\c}}^2 ]\nn\\
&\overset{(b)}{ = }& \E_{i^{t}}[\la \x-\x^{t+1}, (\nabla_{i^t} f(\x^t) + \partial_{i^t} h(\x^{t+1}) -F(\x^t) \cdot  \partial_{i^t} g(\x^{t+1})) \cdot\eit   \ra ] \nn\\
&& -  \E_{i^{t}}[\frac{1}{2} \la \x^{t}- \x^{t+1}, (\bar{\c}_{i^t}  (\x^{t}_{i^t}- \x^{t+1}_{i^t})) \cdot \eit \ra ]\nn\\
&\overset{(c)}{ = }&  \frac{1}{n}\la \x-\x^{t+1}, \nabla f(\x^t) + \partial h(\x^{t+1}) -F(\x^t) \partial g(\x^{t+1}) \ra  -  \frac{1}{2n} \la \x^{t}- \x^{t+1}, \bar{\c}\odot (\x^{t}-\x^{t+1})\ra \nn\\
&\overset{(d)}{ = }&  \frac{1}{n}\la \x-\x^{t+1}, \nabla f(\x^t) + \partial h(\x^{t+1}) \ra\nn\\
 && + \frac{F(\x^t)}{n} [g(\x^{t+1}) - g(\x) + \frac{\bar{\rho}}{2} \| \x-\x^{t+1}\|_{\bar{\c}}^2]    -  \frac{1}{2n}\|\x^{t}- \x^{t+1}\|_{\bar{\c}}^2,
\eeq
\noi where step $(a)$ uses the Pythagoras relation that: $\forall \x,\y, \z,\frac{1}{2}\|\y-\z\|_2^2 - \frac{1}{2}\|\x-\z\|_2^2= \frac{1}{2}\|\x-\y\|_2^2 + \la \y-\x,\x - \z \ra$; step $(b)$ uses the optimality condition in (\ref{eq:pcd:convex:optimality:condition}); step $(c)$ uses the fact that $\E_{i^{t}}[\la \x_{i^t} \eit, \y \ra] = \frac{1}{n} \sum_{j=1}^n  \x_j \y_j = \frac{1}{n}\la \x,\y\ra$; step $(d)$ uses the weakly convexity of $g(\cdot)$ convex that:
\beq  \label{eq:pcd:weak:convex}
-\la  \x- \x^{t+1} ,\partial g(\x^{t+1}) \ra   &\leq & - g(\x) + g(\x^{t+1}) +   \frac{\rho}{2}\|\x-\x^{t+1}\|_2^2. \label{eq:2:weak:p} \nn \\
&  \leq&  - g(\x) + g(\x^{t+1}) +   \frac{\rho}{2 \min(\bar{\c})} \| \x-\x^{t+1}\|_{\bar{\c}}^2\nn \\
&  = &  - g(\x) + g(\x^{t+1}) +   \frac{\bar{\rho}}{2}\| \x-\x^{t+1}\|_{\bar{\c}}^2.\nn
\eeq

We further derive the following results:
\beq \label{eq:parametric:cd:convex:1}
&&\la \x-\x^{t+1}, \nabla f(\x^t) + \partial h(\x^{t+1}) \ra + F(\x^t) (g(\x^{t+1}) - g(\x))\nn\\
&\overset{}{ = } & \la \x-\x^{t+1}, \partial h(\x^{t+1}) \ra + \la \x-\x^{t}, \nabla f(\x^t) \ra + \la \x^{t}-\x^{t+1}, \nabla f(\x^t) \ra + F(\x^t) (g(\x^{t+1}) - g(\x))\nn\\
&\overset{(a)}{ = } & h(\x) - h(\x^{t+1})  + f(\x) - f(\x^{t})  + \la \x^{t}-\x^{t+1}, \nabla f(\x^t) \ra + F(\x^t) (g(\x^{t+1}) - g(\x))\nn\\
&\overset{(b)}{ = } & h(\x) - h(\x^{t+1})  + f(\x)  - f(\x^{t+1}) + \frac{1}{2} \|\x^{t+1} - \x^t\|_{\c}^2 + F(\x^t) (g(\x^{t+1}) - g(\x))\nn\\
&\overset{(c)}{ = } & g(\x) (F(\x) - F(\x^t))  - h(\x^{t+1})     - f(\x^{t+1}) + \frac{1}{2} \|\x^{t+1} - \x^t\|_{\c}^2 + F(\x^t) g(\x^{t+1})\nn\\
&\overset{(d)}{ = } & g(\x) (F(\x) - F(\x^t)) + g(\x^{t+1}) (    F(\x^{t}) -F(\x^{t+1}))   + \frac{1}{2} \|\x^{t+1} - \x^t\|_{\c}^2,
\eeq
\noi where step $(a)$ uses the convexity of $f(\cdot)$ and $h(\cdot)$; step $(b)$ uses the fact that the gradient of $f(\cdot)$ is coordinate-wise Lipschitz continuous that: $\la \x^{t}- \x^{t+1},  \nabla f(\x^t) \ra \leq f(\x^t) - f(\x^{t+1}) + \frac{1}{2} \|\x^{t+1} - \x^t\|_{\c}^2$; step $(c)$ uses the equality that: $f(\x)+h(\x) - F(\x^t)g(\x) = g(\x) (F(\x) - F(\x^t))$; step $(d)$ uses the equality that: $- h(\x^{t+1})     - f(\x^{t+1})+ F(\x^t) g(\x^{t+1}) = g(\x^{t+1}) ( -F(\x^{t+1}) + F(\x^{t}))$.

Combining (\ref{eq:parametric:cd:convex:0}) and (\ref{eq:parametric:cd:convex:1}), we have:
\beq \label{eq:parametric:cd:convex:3}
&& \E_{i^{t}}[\frac{1}{2}\|\x^{t+1}-\x\|_{\bar{\c}}^2 -\E[  \frac{1}{2}\|\x^{t}-\x\|_{\bar{\c}}^2]\nn\\
 &\leq &  \frac{F(\x^t) \bar{\rho}}{2 n}\|\x^{t+1}-\x\|_{\bar{\c}}^2+ g(\x) (F(\x) - F(\x^t)) + g(\x^{t+1}) (    F(\x^{t}) -F(\x^{t+1}))  \nn\\
&\overset{(a)}{ \leq } &  \frac{\varpi}{2 n}\|\x^{t+1}-\x\|_{\bar{\c}}^2+ \frac{g(\x)}{n} (F(\x) - F(\x^t)) + \frac{g(\x^{t+1})}{n} (    F(\x^{t}) -F(\x^{t+1}))  \nn\\
&\overset{(b)}{ \leq } &  \frac{\varpi}{2 n}\|\x^{t+1}-\x\|_{\bar{\c}}^2+ \frac{g(\x)}{n} (F(\x) - F(\x^t)) + \frac{\bar{g}}{n} (    F(\x^{t}) -F(\x^{t+1}))  \nn\\
&\overset{(c)}{ = } &  \frac{\varpi}{2} \E_{i^{t}}[\|\x^{t+1}-\x\|_{\bar{\c}}^2] - \frac{n-1}{n} \frac{\varpi}{2}\|\x^{t}-\x\|_{\bar{\c}}^2 \nn\\
&& + \frac{g(\x)}{n} (F(\x) - F(\x^t)) + \frac{\bar{g}}{n} (    F(\x^{t}) -F(\x^{t+1})),~~~~~~~~
\eeq
\noi where step $(a)$ uses the definition of $\varpi\triangleq F(\x^0) \bar{\rho}$; step $(b)$ uses the assumption that $g(\x^{t})\leq \bar{g},\forall t$; step $(c)$ uses the inequality in (\ref{eq:imp:exp}).

\noi We apply (\ref{eq:parametric:cd:convex:3}) with $\x=\dot{\x}$ and rearranging terms, we obtain:
\beq\label{eq:parametric:cd:convex:34}
\E_{i^{t}}[(1-\varpi) r^{t+1} ] + \frac{\bar{g}}{n} \dot{q}^{t+1} \leq (1-\varpi) r^t + \frac{\varpi}{n}r^t - \frac{g(\x)}{n}\dot{q}^t + \frac{\bar{g}}{n} \dot{q}^t.
\eeq

We now discuss the case when $F(\cdot)$ satisfies the Luo-Tseng error bound assumption. Since $\varpi\leq 1$, we have form (\ref{eq:parametric:cd:convex:34}):
\beq
  \frac{\bar{g}}{n}\dot{q}^{t+1} - \frac{\bar{g}}{n}\dot{q}^{t} + \frac{g(\x)}{n}\dot{q}^t &\leq& (1-\varpi) r^t + \frac{\varpi}{n} r^t    \nn\\
&\overset{(a)}{ \leq } &  (1 + \frac{1}{n}) r^t  \nn\\
&\overset{(b)}{ \leq } &  (1 + \frac{1}{n})  \max(\bar{\c})\delta^2 \frac{\bar{g}}{\theta}(\dot{q}^t-\dot{q}^{t+1})  \nn\\
&\overset{(c)}{ = } &  \kappa_1 \frac{\bar{g}}{n}(\dot{q}^t-\dot{q}^{t+1}) , \label{eq:conv:rate:tsengluo22}
\eeq
\noi where step $(a)$ uses the fact that $0<\varpi\leq 1$; step $(b)$ uses the upper bound for $r^t$ which can be derived using the same strategy as in (\ref{eq:conv:rate:tsengluo:1}); step $(c)$ uses the definition of $\kappa_3$ that $\kappa_1  \triangleq (n + 1) \max(\bar{\c})\delta^2 \frac{1}{\theta}$.

Finally, using the definition of $\kappa_0$ that $\kappa_0 \triangleq \frac{g(\bar{\x})}{\bar{g}}$, we obtain the following results from (\ref{eq:conv:rate:tsengluo22}):
\beq
&& \dot{q}^{t+1} -  \dot{q}^{t} + \kappa_0  \dot{q}^t \leq \kappa_1 (\dot{q}^t-\dot{q}^{t+1}) \nn\\
&\Rightarrow & \dot{q}^{t+1}  \leq \frac{\kappa_1+1 - \kappa_0}{\kappa_1 + 1}\dot{q}^t\nn\\
&\Rightarrow & \dot{q}^{t+1}  \leq (\frac{\kappa_1+1 - \kappa_0}{\kappa_1 + 1})^{t+1}\dot{q}^0.\nn
\eeq

\end{proof}

\section{Proofs for Section \ref{sebsect:concave}} \label{app:sebsect:concave}


\subsection{Proof of Proposition \ref{proposition:global}}
\begin{proof}

\textbf{(a)} We now prove that $F(\cdot)$ is quasi-convex.

First, we prove the following important inequality:
\beq \label{eq:abcd:max}
\frac{a + b}{c+d} \leq \max(\frac{a}{c},\frac{b}{d}),~\forall a\geq 0,b\geq 0,c>0,d>0.
\eeq
\noi We consider two cases. (i) $\frac{a}{c}\leq \frac{b}{d}$. We have $a\leq \frac{b c}{d} \Rightarrow \frac{a+b}{c+d} \leq \frac{\frac{b c}{d}+b}{c+d} = \frac{b}{d} \cdot \frac{c+d}{c+d} = \frac{b}{d} $. (ii) $\frac{a}{c}>\frac{b}{d}$. We have $b < \frac{ad}{c}\Rightarrow \frac{a+b}{c+d} < \frac{a+\frac{ad}{c}}{c+d} = \frac{a}{c} \cdot \frac{c + d}{c+d} = \frac{a}{c}$. Therefore, the inequality in (\ref{eq:abcd:max}) holds.

We derive the following results:
\beq
&& F(\alpha \x + (1-\alpha) \y) \nn\\
&\overset{(a)}{ = }& \frac{f(\alpha \x + (1-\alpha) \y) + h(\alpha \x + (1-\alpha) \y)}{ g(\alpha \x + (1-\alpha) \y)} \nn\\
&\overset{(b)}{ \leq }& \frac{ \alpha f(\x) + (1-\alpha)f(\y) + \alpha h(\x) + (1-\alpha)h(\y) }{ g(\alpha \x + (1-\alpha) \y)}  \nn\\
&\overset{(c)}{ \leq }& \frac{ \alpha ( f(\x) + h(\x)) + (1-\alpha)(f(\y)+h(\y)) }{ \alpha g(\x) + (1-\alpha) g(\y)}   \nn\\
&\overset{(d)}{ \leq }& \max(\frac{ \alpha (f(\x) +h(\x)) }{\alpha g(\x)},\frac{(1-\alpha)(f(\y)  + h(\y))}{(1-\alpha) g(\y)}) \nn\\
&\overset{(e)}{ = }& \max(F(\x),F(\y)), \nn
\eeq
\noi where step $(a)$ uses the definition of $F(\x)$; step $(b)$ uses the convexity of $f(\x)$ and $h(\x)$ that:
\beq
f(\alpha \x + (1-\alpha) \y) \leq \alpha f(\x) + (1-\alpha)f(\y); \nn \\
h(\alpha \x + (1-\alpha) \y) \leq \alpha h(\x) + (1-\alpha)h(\y); \nn
\eeq
\noi step $(c)$ uses the concavity of $g(\x)$ that:
\beq
g(\alpha \x + (1-\alpha) \y) \geq \alpha g(\x) + (1-\alpha)g(\y); \nn
\eeq
\noi step $(d)$ uses the conclusion in (\ref{eq:abcd:max}); step $(e)$ uses the definition of $F(\x)$.

\textbf{(b)} We now prove that any critical point $\bar{\x}$ is also the global optimal solution.

 Assume that $\bar{\x}$ is a critical point of Problem (\ref{eq:main}). We have:
\beq \label{eq:critical:point}
0 \in \nabla f(\bar{\x}) + \partial h(\bar{\x}) - F(\bar{\x}) \partial g(\bar{\x}).
\eeq

Using the convexity of $f(\cdot)$ and $h(\cdot)$, we obtain:
\begin{align}
&~f(\bar{\x}) + h(\bar{\x})\nn \\
\leq &~ f(\x) + h(\x) + \la \bar{\x}-\x, \nabla f(\bar{\x}) + \partial  h(\bar{\x}) \ra  \nn\\
\overset{(a)}{ = } &~ f(\x) + h(\x) + \la \bar{\x}-\x, F(\bar{\x}) \partial g(\bar{\x}) \ra \nn \\
\overset{(b)}{ \leq } & ~f(\x) + h(\x) +  F(\bar{\x}) g(\bar{\x}) - F(\bar{\x}) g(\x) \nn \\
\overset{(c)}{ = } &~ f(\x) + h(\x) +  f(\bar{\x}) + h(\bar{\x}) - F(\bar{\x}) g(\x) , \label{quasiconvex:global}
\end{align}
\noi where step $(a)$ uses the optimality condition in (\ref{eq:critical:point}); step $(b)$ uses the concavity of $g(\cdot)$ that:
\beq
-g(\bar{\x}) \leq -g(\x) + h(\x) - \la \bar{\x}-\x, \partial  g(\bar{\x}) \ra;\nn
\eeq
\noi step $(c)$ uses $F({\x}) g({\x}) = f({\x})+h({\x})$ for all $\x$. Rearranging terms of (\ref{quasiconvex:global}) yields:
\beq
F(\bar{\x}) \leq F(\x),~\forall \x.\nn
\eeq
\noi Thus, we finish the proof of this proposition.

\end{proof}

%
%
%
%
%
%
%

\subsection{Proof of Theorem \ref{the:rate:FCD:PCD:concave}}

\begin{proof}

\textbf{(a)} We prove the convergence rate of \textit{\textbf{FCD}} for convex-concave FMPs.

First, using the first-order optimality condition, we have:
\beq
&&0 \in \frac{    [\nabla_{i^t} f(\x^{t}) + (\c_{i^t} +\theta) \bar{\eta} + \partial_{i^t} h(\x^{t} + \bar{\eta} \ei )  ] - \J_{i^t}(\x^t,\bar{\eta}^t,\theta) \partial_{i^t} g(\x^{t+1})    }{g(\x^{t+1})} \nn\\
&\Leftrightarrow&0 \in      \nabla_{i^t} f(\x^{t}) + (\c_{i^t} +\theta) \bar{\eta} + \partial_{i^t} h(\x^{t+1})    - \alpha^t \partial_{i^t} g(\x^{t+1}) . \label{eq:fractional:cd:concave:optimal:condition}
\eeq

\noi Since $f(\cdot)$ is convex, we have:
\beq
\la \bar{\x}-\x^{t},\nabla f(\x^{t}) \ra \leq f(\bar{\x}) - f(\x^t).\nn
\eeq

Using the fact that $\nabla f(\cdot)$ is coordinate-wise Lipschitz continuous, we have:
\beq
\la \x^t - \x^{t+1},   \nabla f(\x^{t})  \ra   \leq  f(\x^t) - f(\x^{t+1}) + \frac{\c_{i^t}}{2}\|\x^t-\x^{t+1}\|_2^2.\nn
\eeq

Adding the two inequalities above together, we have:
\beq \label{eq:fractional:cd:concave:xt:x:gxt}
 \la \bar{\x}-\x^{t+1},~\nabla f(\x^{t}) \ra  \leq  f(\bar{\x}) - f(\x^{t+1})   +   \frac{\c_{i^t}}{2  }\|\x^t-\x^{t+1}\|_2^2.
\eeq

We derive the following inequalities:
\beq
&& \E_{i^{t}}[\frac{1}{2}\|\x^{t+1}-\x^{t}\|_{\bar{\c}}^2] + \E_{i^{t}}[\frac{1}{2}\|\x^{t+1}-\bar{\x}\|_{\bar{\c}}^2] -  \E[\frac{1}{2}\|\x^{t}-\bar{\x}\|_{\bar{\c}}^2]  \nn\\
&\overset{(a)}{ = }&\E_{i^{t}}[ \la\bar{\x}- \x^{t+1}, \bar{\c}\odot (\x^{t}-\x^{t+1}) \ra ]  \nn\\
&\overset{(b)}{ = }& \E_{\xi^{t}}[\la \bar{\x}-\x^{t+1} , \left( \nabla_{i^t} f(\x^{t}) + \partial_{i^t} h(\x^{t+1})    - \alpha^t \partial_{i^t} g(\x^{t+1}) \right)\eit \ra ] \nn\\
&\overset{(c)}{ = }& \frac{1}{n} \la \bar{\x}-\x^{t+1},~ \nabla f(\x^{t}) + \partial h(\x^{t+1})    - \alpha^t \partial g(\x^{t+1}) \ra \nn\\
&\overset{}{ = }& \frac{1}{n} \la \bar{\x}-\x^{t+1},~  \partial h(\x^{t+1})   \ra +  \frac{1}{n} \la \bar{\x}-\x^{t+1}, \nabla f(\x^{t})  \ra - \frac{ \alpha^t}{n} \la \bar{\x}-\x^{t+1},   \partial g(\x^{t+1})\ra \nn\\
& \overset{(d)}{ \leq } & \frac{1}{n}\left( h(\bar{\x}) - h(\x^{t+1}) \right) +   \frac{1}{n} \la \x-\x^{t+1}, \nabla f(\x^{t})  \ra - \frac{ \alpha^t}{n} \la \bar{\x}-\x^{t+1},   \partial g(\x^{t+1})\ra  \nn\\
& \overset{(e)}{ \leq } & \frac{1}{n}\left( h(\bar{\x}) - h(\x^{t+1}) \right) +   \frac{1}{n} \la \bar{\x}-\x^{t+1}, \nabla f(\x^{t})  \ra + \frac{ \alpha^t}{n}  \left(g(\x^{t+1})- g(\bar{\x})   \right)  \nn\\
& \overset{(f)}{ = } & \frac{1}{n}\left( h(\bar{\x}) - h(\x^{t+1})  +   f(\bar{\x}) - f(\x^{t+1})+\frac{\c_{i^t}}{2}\|\x^t-\x^{t+1}\|_2^2  +  \alpha^t  (g(\x^{t+1})- g(\bar{\x}))   \right),~~~~~\label{eq:fractional:cd:concave:1}
\eeq
\noi where step $(a)$ uses the Pythagoras relation that: $\forall \x,\y,\z, \frac{1}{2}\|\y-\x\|_2^2 + \frac{1}{2}\|\y-\z\|_2^2 - \frac{1}{2}\|\x-\z\|_2^2 = \la \z - \y,\x-\y \ra$; step $(b)$ uses the optimality condition as in (\ref{eq:fractional:cd:concave:optimal:condition}); step $(c)$ uses the fact that $\E_{i^{t}}[\x_{i^t} \eit, \y] = \frac{1}{n}\la \x,\y \ra$; step $(d)$ uses the convexity of $h(\cdot)$ that:
\beq
\la \bar{\x}-\x^{t+1},\partial h(\x^{t+1}) \ra \leq h(\bar{\x}) - h(\x^{t+1});\nn
\eeq
\noi step $(e)$ uses the concavity of $g(\cdot)$ that:
\beq
\la \x^{t+1}-\bar{\x}, \partial g(\x^{t+1})\ra \leq  g(\x^{t+1})- g(\bar{\x}) ;  \nn
\eeq
\noi step $(f)$ uses the inequality in (\ref{eq:fractional:cd:concave:xt:x:gxt}).

From (\ref{eq:fractional:cd:concave:1}) we have the following inequality:
\beq
&&\E_{i^{t}}[\frac{1}{2}\|\x^{t+1}-\bar{\x}\|_{\bar{\c}}^2] - \E[ \frac{1}{2}\|\x^{t}-\bar{\x}\|_{\bar{\c}}^2] \nn\\
&\leq& \frac{1}{n}\left( h(\bar{\x}) - h(\x^{t+1}) \right) +  \frac{1}{n}\left( f(\bar{\x}) - f(\x^{t+1})\right) + \frac{\alpha^t}{n} \left(g(\x^{t+1})- g(\x)   \right) \nn\\
&=& \frac{1}{n}\left(f(\bar{\x}) + h(\bar{\x})- \alpha^t g(\bar{\x}) \right) -  \frac{1}{n}\left( f(\x^{t+1}) + h(\x^{t+1}) - \alpha^t g(\x^{t+1})\right) \nn\\
&\overset{(a)}{ = }& \frac{g(\bar{\x})}{n}\left( F(\bar{\x}) - \alpha^t \right) -  \frac{g(\x^{t+1})}{n}\left( F(\x^{t+1}) - \alpha^t \right) \nn\\
&\overset{(b)}{ = }& \frac{g(\bar{\x})}{n}\left( F(\bar{\x}) - F(\x^{t+1}) \right) +  \frac{\sigma \bar{g}}{n}\left( F(\x^{t}) - F(\x^{t+1}) \right) , \label{eq:fractional:cd:concave:2}
\eeq
\noi where step $(a)$ uses the fact that $F(\x) g(\x) = f(x) + h(\x)$; step $(b)$ uses the Lemma \ref{lemma:omega:sandwich} that: $\alpha^t\geq F(\x^{t+1})$, $\alpha^t - F(\x^{t+1}) \leq \sigma (F(\x^t) - F(\x^{t+1}))$, and the fact that $g(\x^{t+1})\leq \bar{g}$.

From (\ref{eq:fractional:cd:concave:2}), we obtain:
\beq \label{eq:fractional:cd:concave:3}
\E_{i^{t}}[r^{t+1}]  \leq r^{t}  - \frac{g(\bar{\x})}{n} q^{t+1}  + \frac{\sigma \bar{g}}{n} q^t - \frac{\sigma \bar{g}}{n} q^{t+1}.
\eeq

Summing the inequality in (\ref{eq:fractional:cd:concave:3}) over $j=0,1,...,(t-1)$, we have:
\beq
\E_{\xi^{t-1}}[r^t] - r^0 &\leq& - \frac{g(\bar{\x})}{n} \sum_{j=0}^{t-1} q^{j+1} + \frac{\sigma \bar{g}}{n} (q^0 - q^{t})\nn\\
&\overset{(a)}{ \leq } & - \frac{g(\bar{\x})}{n} t q^t + \frac{\sigma \bar{g}}{n} (q^0 + 0),\nn
\eeq
\noi where step $(a)$ uses the fact that $q^j \geq q^t$ for all $j=0,1,...,t$ and $-q^{t}\leq 0$. Finally, combining with that fact that $r^t\geq 0$, we obtain:
\beq
\E_{\xi^{t-1}}[ q^t] \leq  \frac{n(\sigma \bar{g}  q^0 + r^0)}{t g(\bar{\x})}.\nn
\eeq

\textbf{(b)} We prove the convergence rate of \textit{\textbf{PCD}} for convex-concave FMPs.

First, using the first-order optimality condition, we have:
\beq
&&0 \in \frac{    [\nabla_{i^t} f(\x^{t}) + (\c_{i^t} +\theta) \bar{\eta} + \partial_{i^t} h(\x^{t} + \bar{\eta} \ei )  ] - F(\x^t)  \partial_{i^t} g(\x^{t+1})    }{g(\x^{t+1})} \nn\\
&\Leftrightarrow&0 \in      \nabla_{i^t} f(\x^{t}) + (\c_{i^t} +\theta) \bar{\eta} + \partial_{i^t} h(\x^{t+1})    - F(\x^t)  \partial_{i^t} g(\x^{t+1}) . \label{eq:parametric:cd:concave:optimal:condition}
\eeq

Since $f(\cdot)$ and $h(\cdot)$ are convex, we have:
\beq
\la \x-\x^{t},\nabla f(\x^{t}) \ra \leq f(\bar{\x}) - f(\x^t).\nn
\eeq

Using the fact that $\nabla f(\cdot)$ is coordinate-wise Lipschitz continuous, we have:
\beq
\la \x^t - \x^{t+1},   \nabla f(\x^{t})  \ra   \leq  f(\x^t) - f(\x^{t+1}) + \frac{\c_{i^t}}{2}\|\x^t-\x^{t+1}\|_2^2.\nn
\eeq

Adding these two inequalities together, we have:
\beq \label{eq:parametric:cd:concave:xt:x:gxt}
 \la \x-\x^{t+1},~\nabla f(\x^{t}) \ra  \leq  f(\bar{\x}) - f(\x^{t+1})   +   \frac{\c_{i^t}}{2  }\|\x^t-\x^{t+1}\|_2^2.
\eeq

We derive the following inequalities:
\beq
&& \E_{i^{t}}[\frac{1}{2}\|\x^{t+1}-\x^{t}\|_{\bar{\c}}^2 + \frac{1}{2}\|\x^{t+1}-\bar{\x}\|_{\bar{\c}}^2 -  \frac{1}{2}\|\x^{t}-\bar{\x}\|_{\bar{\c}}^2]  \nn\\
&\overset{(a)}{ = }&\E_{i^{t}}[ \la\bar{\x}- \x^{t+1}, \bar{\c}\odot (\x^{t}-\x^{t+1}) \ra ]  \nn\\
&\overset{(b)}{ = }& \E_{i^{t}}[\la \bar{\x}-\x^{t+1} , \left( \nabla_{i^t} f(\x^{t}) + \partial_{i^t} h(\x^{t+1})    - F(\x^t) \partial_{i^t} g(\x^{t+1}) \right)\eit \ra ] \nn\\
&\overset{(c)}{ = }& \frac{1}{n} \la \bar{\x}-\x^{t+1},~ \nabla f(\x^{t}) + \partial h(\x^{t+1})    - F(\x^t) \partial g(\x^{t+1}) \ra \nn\\
&\overset{}{ = }& \frac{1}{n} \la \bar{\x}-\x^{t+1},~  \partial h(\x^{t+1})   \ra +  \frac{1}{n} \la \bar{\x}-\x^{t+1}, \nabla f(\x^{t})  \ra - \frac{ F(\x^t)}{n} \la \x-\x^{t+1},   \partial g(\x^{t+1})\ra \nn\\
& \overset{(d)}{ \leq } & \frac{1}{n}\left( h(\bar{\x}) - h(\x^{t+1}) \right) +   \frac{1}{n} \la \bar{\x}-\x^{t+1}, \nabla f(\x^{t})  \ra - \frac{ F(\x^t)}{n} \la \bar{\x}-\x^{t+1},   \partial g(\x^{t+1})\ra  \nn\\
& \overset{(e)}{ \leq } & \frac{1}{n}\left( h(\bar{\x}) - h(\x^{t+1}) \right) +   \frac{1}{n} \la \bar{\x}-\x^{t+1}, \nabla f(\x^{t})  \ra + \frac{ F(\x^t)}{n}  \left(g(\x^{t+1})- g(\bar{\x})   \right)  \nn\\
& \overset{(f)}{ = } & \frac{1}{n}\left( h(\bar{\x}) - h(\x^{t+1}) +   f(\bar{\x}) - f(\x^{t+1})+\frac{\c_{i^t}}{2}\|\x^t-\x^{t+1}\|_2^2 +  F(\x^t) (g(\x^{t+1})- g(\bar{\x})   ) \right),~~~~~~~~ \label{eq:parametric:cd:concave:0}
\eeq
\noi where step $(a)$ uses the Pythagoras relation that: $\forall \x,\y,\z, \frac{1}{2}\|\y-\x\|_2^2 + \frac{1}{2}\|\y-\z\|_2^2 - \frac{1}{2}\|\x-\z\|_2^2 = \la \z - \y,\x-\y \ra$; step $(b)$ uses the optimality condition as in (\ref{eq:parametric:cd:concave:optimal:condition}); step $(c)$ uses the fact that $\E_{i^{t}}[\x_{i^t} \eit, \y] = \frac{1}{n}\la \x,\y \ra$; step $(d)$ uses the convexity of $h(\cdot)$ that:
\beq
\la \x-\x^{t+1},\partial h(\x^{t+1}) \ra \leq h(\x) - h(\x^{t+1});\nn
\eeq
\noi step $(e)$ uses the concavity of $g(\cdot)$ that:
\beq
\la \x^{t+1}-\x, \partial g(\x^{t+1})\ra \leq  g(\x^{t+1})- g(\x)  ; \nn
\eeq
\noi step $(f)$ uses the inequality in (\ref{eq:parametric:cd:concave:xt:x:gxt}).

We have the following inequalities:
\beq
 f({\x}) + h({\x}) - F(\x^t) g({\x}) \leq g({\x}) ( \frac{f({\x})+g({\x} )}{g({\x})} - F(\x^t) ) \leq g({\x}) ( F({\x})  - F(\x^{t}) )\label{eq:parametric:cd:concave:fhg:x}\\
-h(\x^{t+1}) -f(\x^{t+1}) + F(\x^{t}) g(\x^{t+1}) =   g(\x^{t+1}) (F(\x^{t}) - F(\x^{t+1})) . \label{eq:parametric:cd:concave:fhg:xt1}
\eeq

Combining (\ref{eq:parametric:cd:concave:0}), (\ref{eq:parametric:cd:concave:fhg:x}), and (\ref{eq:parametric:cd:concave:fhg:xt1}), we obtain:
\beq
\E_{i^{t}}[ \frac{1}{2}\|\x^{t+1}-\x\|_{\bar{\c}}^2]      \leq  \frac{1}{2}\|\x^{t}-\x\|_{\bar{\c}}^2 + \frac{g(\bar{\x})}{n} ( F(\bar{\x})  - F(\x^{t}) )  + \frac{g(\x^{t+1})}{n} (F(\x^{t}) - F(\x^{t+1})) . \label{eq:parametric:cd:concave:2}
\eeq

Using (\ref{eq:parametric:cd:concave:2}) and the fact that $\g(\x^t)\leq \bar{g}$, we obtain:
\beq \label{eq:parametric:cd:concave:3}
\E_{i^{t}}[r^{t+1}]  \leq r^{t}  - \frac{g(\bar{\x})}{n} q^{t}  + \frac{ \bar{g}}{n} q^t - \frac{ \bar{g}}{n} q^{t+1}.\nn
\eeq

Summing the inequality above over $j=0,1,...,t$, we have:
\beq
\E_{\xi^{t}}[r^{t+1}] - r^{0}  &\leq& - \frac{g(\bar{\x})}{n} \sum_{j=0}^{t} q^j + \frac{\bar{g}}{n} (q^0 - q^{t+1}) \nn\\
& \overset{(a)}{ \leq } & \frac{g({\bar{\x}})}{n} - \sum_{j=0}^{t} q^t + \frac{\bar{g}}{n} q^0 \nn\\
& \overset{}{ = } & \frac{- g({\bar{\x}})}{n} (t+1)q^t + \frac{\bar{g}}{n}q^0 ,\nn
\eeq
\noi where step $(a)$ uses $q^j \geq q^t$ for all $j=0,1,...,t$ and $-q^{t+1}\leq 0$. Finally, we have the following result:
\beq
\E_{\xi^{t-1}}[q^t] \leq  \frac{ \bar{g} nq^0 + nr^0}{g(\bar{\x}) (t+1) }.\nn
\eeq
\end{proof}

\section{Additional Discussions} \label{app:sebsect:discussion}

In this section, we discuss the optimality hierarchy, the globally/locally bounded non-convexity assumption, and the convexity of the function $g(\x)=\|\G\x\|_4^2$.

\subsection{Fractional Reformulations for Problem (\ref{eq:pca0})}
\label{app:sect:disc:ref}

First, we focus on the following minimization problems with $\QQ\succ \mathbf{0}$:
\beq
&& \bar{\v} = \arg \min_{\v} ~F_1(\v) \triangleq - \|\G\v\|_p,~s.t.~\v^T\QQ\v = 1 \label{eq:000} \\
&& \bar{\x} = \arg  \min_{\x}~ F_2(\x) \triangleq \frac{\x^T\QQ\x +\gamma_1}{\|\G\x\|_p + \gamma_2}\label{eq:111}
.
\eeq

The following proposition establish the relations between Problem (\ref{eq:000}) and Problem (\ref{eq:111}).

\begin{proposition}  \label{prop:equ:lppp}

We have the following results.

\noi \textbf{(a)} If $\bar{\v}$ is an optimal solution to $(\ref{eq:000})$, then ${ \pm \bar{\alpha} \bar{\v}}$ with $\bar{\alpha}  \in \arg \min_{\alpha}~ \frac{\bar{\v}^T\QQ\bar{\v} \alpha^2 +\gamma_1}{\alpha \|\G(\bar{\v})\|_p + \gamma_2}$ is an optimal solution to $(\ref{eq:111})$.

\noi\textbf{(b)} If $\bar{\x}$ is an optimal solution to $(\ref{eq:111})$, then ${ \pm \bar{\beta} \bar{\v}}$ with $\bar{\beta} = \pm 1 / \sqrt{\bar{\x}^T\QQ\bar{\x}}$ is an optimal solution to $(\ref{eq:000})$.

\begin{proof}

We notice that Problem (\ref{eq:000}) can be rewritten as:
\beq
\bar{\v} = \arg \min_{\v} ~F_1(\v) \triangleq - \|\G\v\|_p,~s.t.~\v^T\QQ\v \leq 1 \label{eq:000:b}
\eeq

\noi On the one hand, since $\bar{\v}$ is an optimal solution to Problem (\ref{eq:000:b}), there exists a multiplier $\theta_1>0$ which is associated to the constraint $\v^T\QQ\v \leq 1$ as in Problem (\ref{eq:000:b}) that:
\beq\label{eq:000:z}
\bar{\v} = \arg \min_{\v}~ \acute{F} (\v) \triangleq \theta_1(\v^T\QQ\v-1) - \|\G\v\|_p.\nn
\eeq
\noi On the other hand, since $\bar{\x}$ is an optimal solution to Problem (\ref{eq:111}), there exists a constant $\theta_2>0$ that:
\beq\label{eq:111:z}
\bar{\x} = \arg  \min_{\x}~ \check{F}(\x) \triangleq (\x^T\QQ\x +\gamma_1) - \theta_2(\|\G\x\|_p + \gamma_2)\nn
\eeq

\noi It is not hard to notice that the gradient of $\acute{F} (\v)$ and $\check{F}(\x)$ can be computed as:
\begin{align}
\nabla \acute{F} (\v) &~= 2 \theta_1 \QQ\v - \|\G\v\|_p^{1-p} \G^T ( \text{sign}(\G\v) \odot |\G\v|^{p-1}). \nn \\
\nabla \acute{F} (\x) &~= 2 \QQ\x -  \theta_2 \|\G\x\|_p^{1-p} \G^T ( \text{sign}(\G\x) \odot |\G\x|^{p-1}). \nn
\end{align}
\noi By the first-order optimality condition, we have:
\begin{align}
\v &~= \frac{1}{2\theta_1}\QQ^{-1}\left(\|\G\v\|_p^{1-p} \G^T ( \text{sign}(\G\v) \odot |\G\v|^{p-1})\right) , \label{eq:000:zz} \\
\x &~= \frac{\theta_2}{2}\QQ^{-1}\left(\|\G\x\|_p^{1-p} \G^T ( \text{sign}(\G\x) \odot |\G\x|^{p-1})\right) .\label{eq:111:zz}
\end{align}
\noi In view of (\ref{eq:000:zz}) and (\ref{eq:111:zz}), we conclude that the optimal solution for Problem (\ref{eq:000}) and Problem (\ref{eq:111}) only differ by a scale factor.

\noi \textbf{(a)} Since $\bar{\v}$ is the optimal solution to $(\ref{eq:000})$, the optimal solution to Problem $(\ref{eq:111})$ can be computed as $\bar{\alpha} \cdot \bar{\v}$ with
\begin{align}
\bar{\alpha} &~= \arg \min_{\alpha}~F_2(\alpha \cdot \bar{\v}) \nn\\
&~=\arg \min_{\alpha}~ \frac{\bar{\v}^T\QQ\bar{\v} \alpha^2 +\gamma_1}{\alpha \|\G(\bar{\v})\|_p + \gamma_2}. \nn
\end{align}

\noi \textbf{(b)} Since $\bar{\x}$ is the optimal solution to $(\ref{eq:111})$, the optimal solution to Problem $(\ref{eq:000})$ can be computed as $\bar{\beta} \cdot \bar{\x}$ with
\beq
\bar{\beta} &=& \arg \min_{\beta}~F_1(\beta \cdot \bar{\x}) , ~s.t.~(\beta \cdot \bar{\x})^T \QQ (\beta \cdot \bar{\x}) = 1\nn
\eeq
\noi After some preliminary calculations, we have: $\bar{\beta} = \pm 1 / \sqrt{\bar{\x}^T\QQ\bar{\x}}$.

\end{proof}

\end{proposition}

Second, we focus on the following minimization problems with $\QQ\succ \mathbf{0}$:
\beq
&& \bar{\v} = \arg \min_{\v} ~F'_1(\v) \triangleq - \|\G\v\|_p^2,~s.t.~\v^T\QQ\v = 1 \label{eq:000:hhh} \\
&& \bar{\x} = \arg  \min_{\x}~ F'_2(\x) \triangleq \frac{\x^T\QQ\x +\gamma_3}{\|\G\x\|^2_p + \gamma_4}\label{eq:111:hhh}.
\eeq
\noi Note that Problem (\ref{eq:000:hhh}) is equivalent to Problem (\ref{eq:000}).

The following proposition establish the relations between Problem (\ref{eq:000:hhh}) and Problem (\ref{eq:111:hhh}).

\begin{proposition}  \label{prop:equ}

We have the following results.

\noi \textbf{(a)} If $\bar{\v}$ is an optimal solution to $(\ref{eq:000:hhh})$, then ${ \pm \bar{\alpha} \bar{\v}}$ with $\bar{\alpha}  \in \arg \min_{\alpha}~ \frac{\bar{\v}^T\QQ\bar{\v} \alpha^2 +\gamma_3}{\|\G(\bar{\v})\|^2_p \alpha^2 + \gamma_4}$ is an optimal solution to $(\ref{eq:111:hhh})$.

\noi\textbf{(b)} If $\bar{\x}$ is an optimal solution to $(\ref{eq:111:hhh})$, then ${ \pm \bar{\beta} \bar{\v}}$ with $\bar{\beta} = \pm 1 / \sqrt{\bar{\x}^T\QQ\bar{\x}}$ is an optimal solution to $(\ref{eq:000:hhh})$.

\begin{proof}

The proof of this proposition is analogous to that of Proposition \ref{prop:equ:lppp}. We omit the proof for brevity.

\end{proof}

\end{proposition}

\begin{figure} [!h]
  \centering
  \includegraphics[width=3.6cm,height=2.7cm]{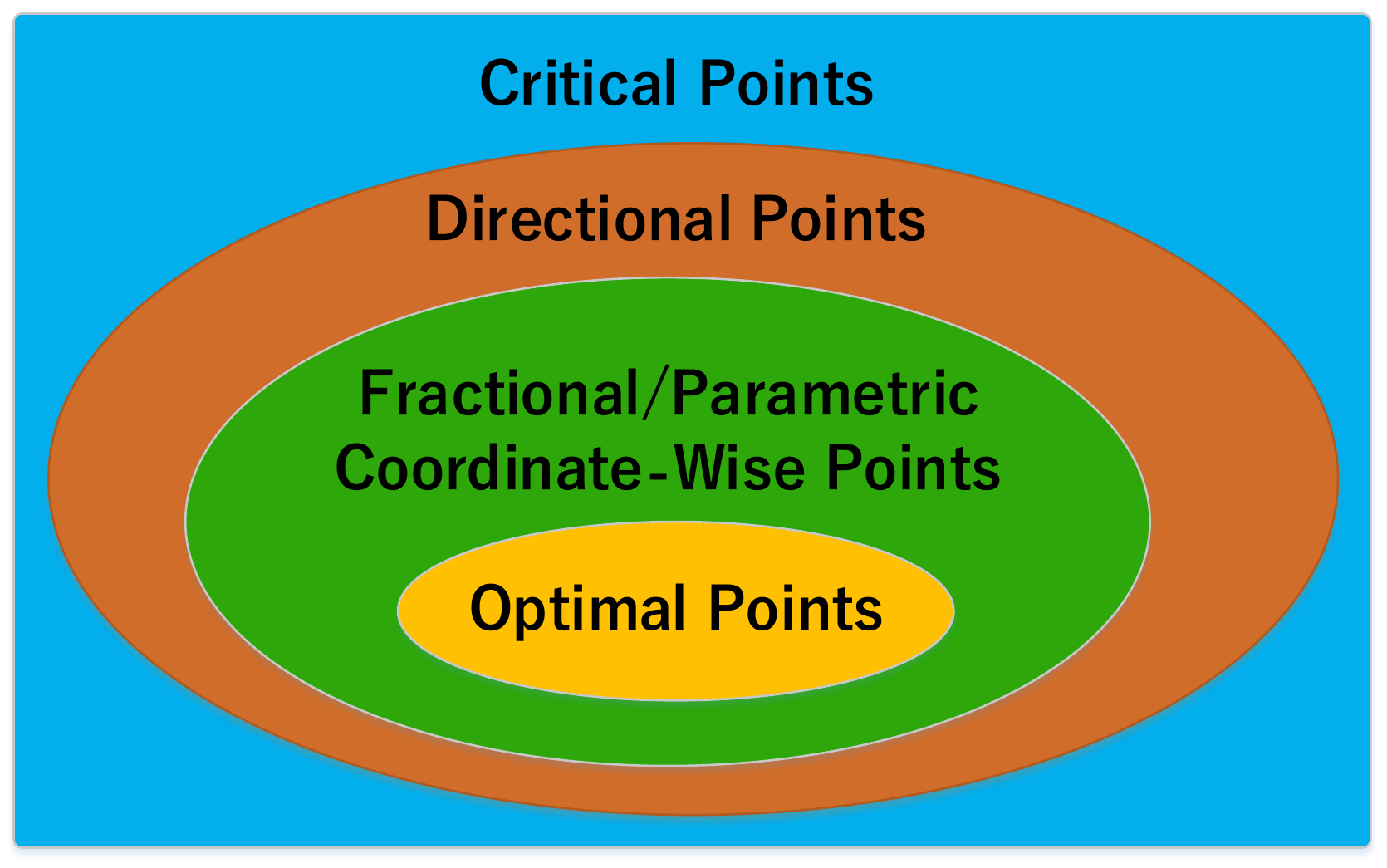}~
  \includegraphics[width=4.6cm,height=2.8cm]{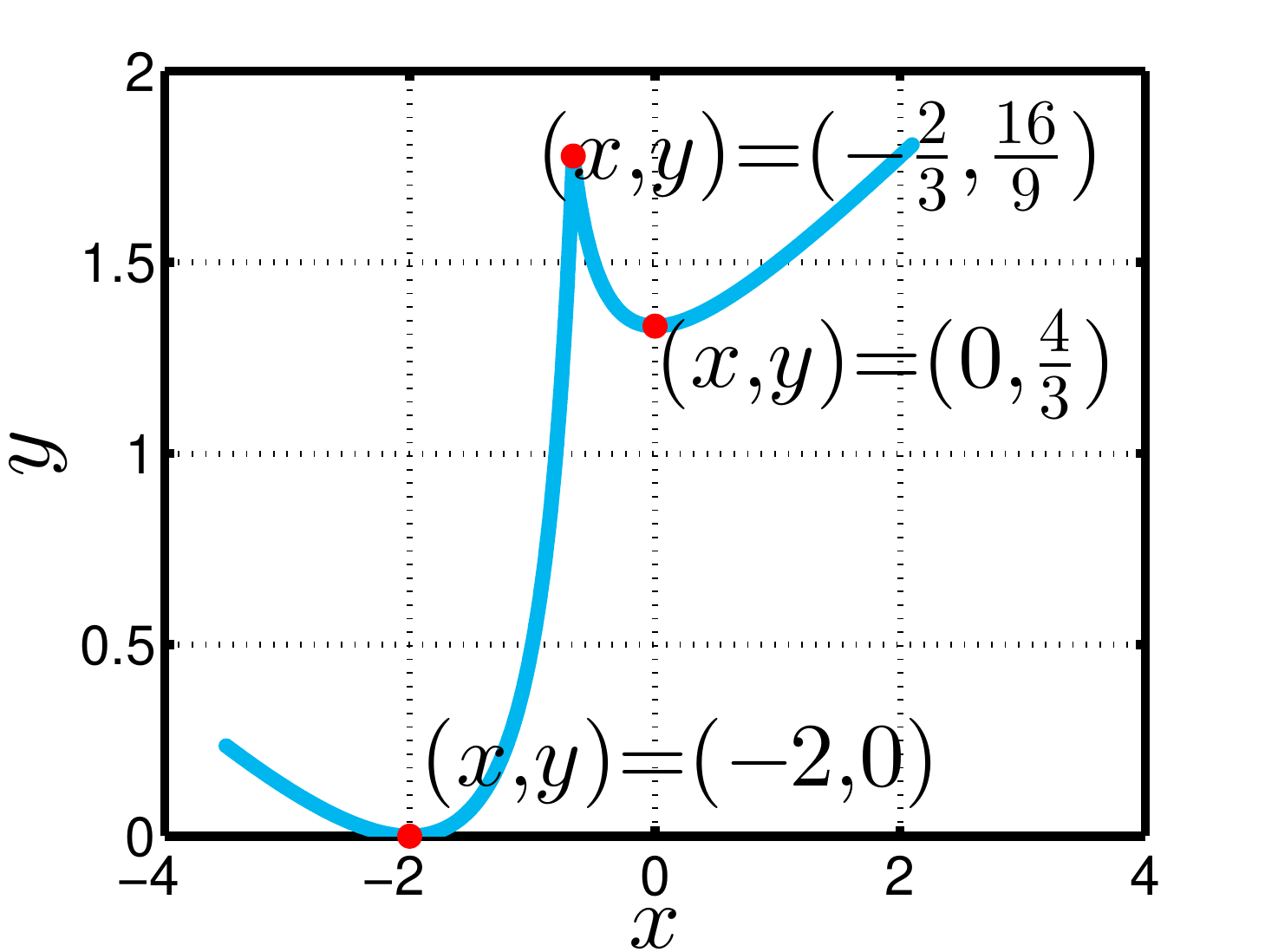}

  \caption{Left: optimality hierarchy between the optimality conditions. Note that the condition of \textit{FCW}-point is equivalent to that of \textit{PCW}-point. Right: Geometric interpretation for the one-dimensional fractional problem with $F(x) \triangleq \frac{(x+2)^2}{|3x+2|+1}$. }
\label{fig:optimal:hierarchy}
\end{figure}

\begin{table}[!h]
\begin{center}
\scalebox{0.7}{\begin{tabular}{|c|c|c|c|c|c|c|}
\hline
$x$    & $F(x)$ & C-point & D-point & \textit{FCW}-point& \textit{PCW}-point \\
\hline
$x_1=-\tfrac{2}{3}$  & $(\tfrac{4}{3})^2$  & \red{Yes}  &No &No&No\\
\hline
$x_2=0 $  & $\frac{4}{3}$  & \red{Yes} & \red{Yes} &No&No\\
\hline
$x_3=-2 $  & 0  & \red{Yes} & \red{Yes} & \red{Yes} &\red{Yes} \\
  \hline
\end{tabular}}
\end{center}
\caption{Points satisfying different optimality conditions.}
\label{table:optimality}
\end{table}
\subsection{A Simple Example for the Optimality Hierarchy}\label{app:sebsect:hierarchy}

To show the optimality hierarchy between the optimality conditions, we consider the following one-dimensional example which has been mentioned in the paper:
\beq
\min_{x} F(\x)  \triangleq \frac{(x+2)^2}{|3x+2|+1}\nn
\eeq

This problem contains three \textit{C}-points $\{-\tfrac{2}{3},0,-2\}$, and we now show that this problem contains one unique Parametric Coordinate-Wise Point (\textit{PCW}-point). \bfit{(i)} We consider the point $x_1=-\frac{2}{3}$. We have the following parametric problem:
\begin{align}
  &~\arg \min_y P_1 (y) \triangleq  (y+2)^2 - F(x_1) (|3y+2|+1) \nn\\
\overset{(a)}{ = }  &~ \arg \min_y P_1 (y) \triangleq (y+2)^2 - (\tfrac{4}{3})^2 (|3y+2|+1) \nn\\
\overset{(b)}{ = } &~ \arg \min_{y} P_1 (y) \triangleq (y+2)^2 - (\tfrac{4}{3})^2 (|3y+2|+1), s.t.~y \in \{-\tfrac{2}{3},0,-4\} \nn\\
\overset{(c)}{ = } &~ -4 \nn  \neq   x_1, \nn
\end{align}
\noi where step $(a)$ uses $F(x_1) = (\frac{4}{3})^2$; step $(b)$ uses the fact that $\{-\tfrac{2}{3},0,-4\}$ are the three critical points of $\min_{y} P_1 (y)$; step $(c)$ uses the fact that $P_1(-\tfrac{2}{3}) = \frac{4}{9}$, $P_1(0)=0$, $P_1(-4) = -\frac{32}{3}$, and $y=-4$ is the global minimizer of the problem $\min_{y} P_1 (y)$. Since $-4 \nn  \neq   x_1 = -\frac{2}{3}$,  $x_1 = -\frac{2}{3}$ is not a \textit{PCW}-point.

\bfit{(ii)} We consider the point $x_2=0$. We have the following parametric problem:
\begin{align}
  &~\arg \min_y P_2 (y) \triangleq  (y+2)^2 - F(x_2) (|3y+2|+1) \nn\\
\overset{(a)}{ = }  &~ \arg \min_y P_1 (y) \triangleq (y+2)^2 - \tfrac{4}{3}(|3y+2|+1) \nn\\
\overset{(b)}{ = } &~ \arg \min_{y} P_1 (y) \triangleq (y+2)^2 - \tfrac{4}{3} (|3y+2|+1), s.t.~y \in \{-\tfrac{2}{3},\tfrac{2}{3},-\tfrac{14}{3}\} \nn\\
\overset{(c)}{ = } &~ -\tfrac{14}{3} \nn  \neq   x_2, \nn
\end{align}
\noi where step $(a)$ uses $F(x_1) = (\frac{4}{3})^2$; step $(b)$ uses the fact that $\{-\tfrac{2}{3},\tfrac{2}{3},-\tfrac{14}{3}\}$ are the three critical points of $\min_{y} P_1 (y)$; step $(c)$ uses the fact that $P_1(-\tfrac{2}{3}) =0$, $P_1(0)=\frac{16}{27}$, $P_1(-4) = -16$, and $y=-\tfrac{14}{3}$ is the global minimizer of the problem $\min_{y} P_2 (y)$. Since $-\tfrac{14}{3}  \nn  \neq   x_2 = 0$,  $x_2 = 0$ is not a \textit{PCW}-point.

\bfit{(iii)} We consider the point $x_3=-2$. We have the following parametric problem:
\begin{align}
  &~\arg \min_y P_3 (y) \triangleq  (y+2)^2 - F(x_3) (|3y+2|+1) \nn\\
\overset{(a)}{ = }  &~ \arg \min_y P_1 (y) \triangleq (y+2)^2 \nn\\
\overset{(b)}{ = } &~ -2 =   x_3 ,\nn
\end{align}
\noi where step $(a)$ uses $F(x_3) = 0$; step $(b)$ uses the fact that $y=-2$ is the global minimizer of the problem $\arg \min_y P_1 (y)$. Since $-2 = x_3$,  $x_3 = -2$ is a \textit{PCW}-point.

Therefore, $x=-2$ is the unique \textit{PCW}-point.

\noi Figure \ref{fig:optimal:hierarchy} demonstrates the optimality hierarchy and the geometric interpretation for the one-dimensional problem. Table \ref{table:optimality} shows the points satisfying different optimality conditions. We conclude that the condition of \textit{FCW}-point and \textit{PCW}-point might be a much stronger condition than the condition of critical point and direction point.

\subsection{The Globally or Locally Bounded Non-Convexity Assumption} \label{ass:globally:locally}

We prove that $\tilde{g}(\x) =  - \|\G\x\|_4^2$ is globally $\rho$-bounded non-convex, while $\tilde{g}(\x) = -\sum_{j=1}^k |\x_{[j]}| $ is locally $\rho$-bounded non-convex.

\begin{lemma}
Assume $\x \neq \mathbf{0}$ and $\A$ has full column rank. The function $\tilde{g}(\x) =  - \|\G\x\|_4^2$ is globally $\rho$-bounded non-convex with $\rho=6m\max_{i} (\G\G^T)_{ii} \cdot \frac{\lambda_{\max}(\G^T\G)}{\lambda_{\min}(\G^T\G)}  $, where $\lambda_{\min}(\G^T\G)$ and $\lambda_{\max}(\G^T\G)>0$ denote the smallest and the largest eigenvalue of the matrix $\G^T\G$, respectively.

\begin{proof}
The first-order and second-order gradient of $\tilde{g}(\x)$ can be computed as:
\beq
\nabla \tilde{g} (\x) &=& \frac{ 2\sum_{i}^m(\G_i\x)^3\G_i^T}{- \|\G\x\|_4^2}, \nn\\
\nabla^2 \tilde{g}(\x)&=&\frac{ 6 \sum_{i}^m [ (\G_i \x)^2 \G_i^T\G_i] \|\G\x\|_4^2}{- \|\G\x\|_4^4} + \frac{ 2\sum_{i}^m(\G_i\x)^3\G_i^T \nabla \tilde{g}(\x)^T }{-\|\G\x\|_4^4}, \nn
\eeq
\noi where $\G\in \mathbb{R}^{m\times n}$ and $\G_i \in \mathbb{R}^{1\times n}$ is the $i$-the row of $\G$.

The $\rho$-bounded nonconvexity of $\tilde{g}(\x)$ is equivalent to the convexity of $(\tilde{g}(\x)+\tfrac{\rho}{2}\|x\|_2^2)$. In what follows, we prove that $\nabla^2 \tilde{g}(\x) + \rho \I\succcurlyeq  \textbf{0}$.

\textbf{(a)} We bound the term $\sum_{i}^m (\G_i \x)^2 \G_i^T\G_i$. We denote $\v_i \triangleq \|\G_i\|_2^2$ with $i=1,2,...,m$. We have the following upper-bound:
\beq
\sum_{i}^m (\G_i \x)^2 \G_i^T\G_i \preceq \sum_{i}^m \|\x\|_2^2 \|\G_i\|_2^2  \G_i^T\G_i = \G^Tdiag(\v)\G \|\x\|_2^2 \preceq \|\G\|_2^2 \max(v)\|\x\|_2^2,\label{eq:GxGx2}
\eeq
\noi where the first inequality uses the Cauchy-Schwarz inequality and the last inequality uses the norm inequality.

\textbf{(b)} We bound the term $\|\G\x\|_4^2$. Using the fact that $\sqrt{m}\|\y\|_4\geq \|\y\|_2\geq \|\y\|_4$ for all $\y \in \mathbf{R}^m$. We have the following lower-bound:
\beq
\|\G\x\|_4^2 \geq \frac{1}{m}\|\G\x\|_2^2 \geq \frac{1}{m}\lambda_{\min}(\G^T\G)\|\x\|_2^2. \label{eq:GxGx3}
\eeq

\textbf{(c)} Finally, we have the following inequalities:
\beq \label{eq:GxGx1}
\nabla^2 \tilde{g}(\x)  &\overset{(a)}{ \succcurlyeq }&   \frac{6\sum_{i}^m [ (\G_i \x)^2 \G_i^T\G_i  ] }{-\|\G\x\|_4^2} + \mathbf{0} \nn\\
 &\overset{(b)}{ \succcurlyeq }&   -6\frac{\|\G\|_2^2 \max(v)\|\x\|_2^2}{\frac{1}{m}\lambda_{\min}(\G^T\G)\|\x\|_2^2} \cdot \I, \nn\\
 &\overset{(c)}{ = }& -6m \max(v) \cdot \frac{\lambda_{\max}(\G^T\G) }{\lambda_{\min}(\G^T\G)} \cdot \I = - \rho \I,\nn
\eeq
\noi where step $(a)$ uses the fact that $\frac{ 2\sum_{i}^m(\G_i\x)^3\G_i^T \nabla \tilde{g}(\x)^T }{-\|\G\x\|_4^4}$ is positive semidefinite, step $(b)$ uses (\ref{eq:GxGx2}) and (\ref{eq:GxGx3}); step $(c)$ uses the definition of $\rho$.

\end{proof}
\end{lemma}

Note that the assumption $\x \neq \mathbf{0}$ is automatically satisfied by Problem (\ref{algo:main}) since we assume that $g(\x)>0$.

\begin{lemma}

The function $\tilde{g}(\x) = -\sum_{j=1}^k |\x_{[j]}| $ is locally $\rho$-weakly convex with $\rho<+\infty$.

\begin{proof}

For simplicity, we define $\|\x\|_{[k]}\triangleq \sum_{j=1}^k |\x_{[j]}|$. For any $\x\in\mathbb{R}^n$ and a given parameter $k$, the subgradient of $\|\x\|_{[k]}$ can be computed as $\partial \|\x\|_{[k]} = {\tiny \{
                               \begin{array}{ll}
                                 \text{sign}(\x_i),  & \hbox{ $i \in \Delta_k(\x)\text{~and~} \x_i \neq 0$;} \\
                                 $[-1,1]$, & \hbox{else.}
                               \end{array}
                             \}}$, where $\Delta_k(\x)$ is the index of the largest (in magnitude) $k$ elements of $\x$.

As the two reference points $\x\neq \y$ in Assumption \ref{ass:1}, we assume that there exists a constant $\epsilon>0$ satisfying $\|\x-\y\|_2 \geq \epsilon $. We have:
\beq
&& \tilde{g}(\x) - \tilde{g}(\y) - \la \x-\y, \partial \tilde{g}(\x)\ra \nn\\
 &= & - \|\x\|_{[k]} + \|\y\|_{[k]} - \la \x-\y,   \partial (-\|\x\|_{[k]})\ra   \nn\\
 &\overset{(a)}{ \leq }&    \|\y-\x\|_{[k]} + \|\y-\x\| \cdot \|\partial (\|\x\|_{[k]})\|   \nn\\
 &\overset{(b)}{ \leq }&   \|\y-\x\|_{1} + \|\y-\x\| \cdot \sqrt{n}\nn\\
 &\overset{(c)}{ \leq }& 2 \sqrt{n}\|\x-\y\|_2 \nn\\
 &\overset{(d)}{ \leq }& \tfrac{2\sqrt{n}}{\epsilon}\|\x-\y\|_2^2, \nn
\eeq
\noi where step $(a)$ uses the triangle inequality that $\|\y\|_{[k]}-\|\x\|_{[k]}\leq \|\y-\x\|_{[k]}$ since $\|\cdot\|_{[k]}$ is a norm; step $(b)$ uses the fact that $ \|\partial (\|\x\|_{[k]})\|\leq \sqrt{n}$; step $(c)$ uses the fact that $\|\x\|\leq\|\x\|_1\leq \sqrt{n}\|\x\|$ for all $\x\in\mathbb{R}^n$; step $(d)$ uses $\|\x-\y\|_2 \geq \epsilon$.

\noi Therefore, the function $\tilde{g}(\x)$ is $\rho$-bounded non-convex with $\rho<+\infty$.

\end{proof}
\end{lemma}


\subsection{The function $g(\x) = \|\G\x\|_4^2$ is convex}

We prove that the function $g(\x) = \|\G\x\|_4^2$ is convex. We first present the following useful lemma.

\begin{lemma} \label{lemma:nonnegative:convex}
Assume that $p(\x)$ is a convex and non-negative function. The function $g(\x) = (p(\x))^2$ is convex.
\begin{proof}

By the convexity of $p(\x)$, we have:
\beq
p((1-t)\x+t\y) \leq (1-t)p(\x)+tp(y),\forall t \in (0,1). \nn
\eeq
\noi Squaring both sides of the inequality above, we obtain:
\beq
&& p((1-t)\x+t\y) \cdot p((1-t)\x+t\y) \nn\\
&\leq& (1-t)(1-t) p(\x) p(\x) + t^2 p(y) p(y) + 2 t(1-t) p(\x)p(\y) \nn\\
&=& (1-t) p(\x) p(\x) - (1-t) t p(\x) p(\x) + t^2 p(y) p(y) + 2 t(1-t) p(\x)p(\y) \nn\\
&=& (1-t) p(\x) p(\x) + t p(y) p(y)  - t(1-t) (p(y)- p(\x))^2 \nn\\
& \leq & (1-t)p(\x)p(\x) + t p(y) p(y)  . \nn
\eeq
\noi where the last step uses $t \in(0,1)$ and $(p(y)- p(\x))^2\geq0$.

\end{proof}
\end{lemma}

Note that $p(\x) = \|\G\x\|_4$ is a convex and non-negative function. Using Lemma \ref{lemma:nonnegative:convex}, we conclude that $g(\x)=\|\G\x\|_4^2$ is convex with respect to $\x$.

\end{document}